\documentclass{birkmult}
\usepackage{color}

 \newtheorem{theorem}{Theorem}[section]
 \newtheorem{corollary}[theorem]{Corollary}
 \newtheorem{lemma}[theorem]{Lemma}
 \newtheorem{proposition}[theorem]{Proposition}
 \theoremstyle{definition}
 \newtheorem{definition}[theorem]{Definition}
 \theoremstyle{remark}

 \numberwithin{equation}{section}
 \newcommand{\cir}{\mbox{$\bigcirc\kern -8pt c\ $}}
 \newcommand{\sph}{\mbox{$\bigcirc\kern -8pt s\ $}}
 \DeclareMathOperator{\im}{Im}
\DeclareMathOperator{\Li}{Li}

\begin{document}
%
%
%
%
%
%
%
%
%
\title[Surfaces from Circles]
 {Surfaces from Circles}
\author[Alexander I. Bobenko]{Alexander I. Bobenko}

\address{%
Institut f\"ur Mathematik,\\
Technische Universit\"at Berlin,\\
Strasse des 17. Juni 136, 10623 Berlin,\\
Germany}

\email{bobenko@math.tu-berlin.de}

\thanks{Partially supported by the DFG Research Unit 565 ``Polyhedral Surfaces''.}




\begin{abstract}
In the search for appropriate discretizations of surface theory it
is crucial to preserve such fundamental properties of surfaces as
their invariance with respect to transformation groups. We discuss
discretizations based on M\"obius invariant building blocks such
as circles and spheres. Concrete problems considered in these
lectures include the Willmore energy as well as conformal and
curvature line parametrizations of surfaces. In particular we
discuss geometric properties of a recently found discrete Willmore
energy. The convergence to the smooth Willmore functional is shown
for special refinements of triangulations originating from a
curvature line parametrization of a surface. Further we treat
special classes of discrete surfaces such as isothermic and
minimal. The construction of these surfaces is based on the theory
of circle patterns, in particular on their variational
description.
\end{abstract}

\maketitle
\section{Why from Circles?}
\label{s.intro}

The theory of polyhedral surfaces aims to develop discrete
equivalents of the geometric notions and methods of smooth surface
theory. The latter appears then as a limit of refinements of the
discretization. Current interest in this field derives not only
from its importance in pure mathematics but also from its
relevance for other fields like computer graphics.

\begin{figure}[h]
\begin{center}
\includegraphics[width=0.50\textwidth]{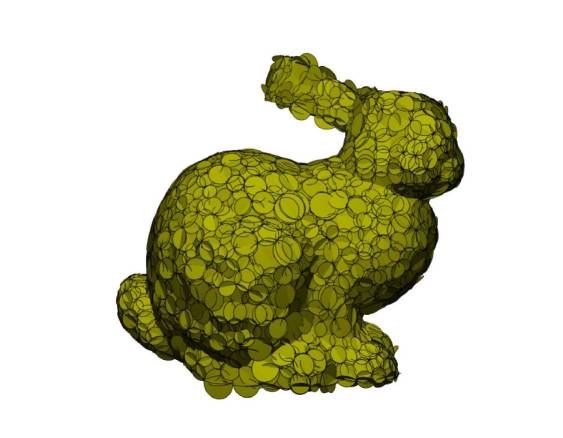}
\includegraphics[width=0.49\textwidth]{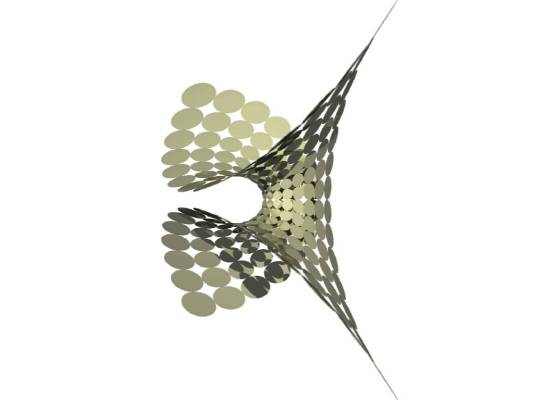}
\end{center}
\caption{Discrete surfaces made from circles: general simplicial surface and a
discrete minimal Enneper surface.
 \label{f.dSurfaces}}
\end{figure}

One may suggest many different reasonable discretizations with the
same smooth limit. Which one is the best? In the search for
appropriate discretizations, it is crucial to preserve the
fundamental properties of surfaces. A natural mathematical
discretization principle is the invariance with respect to
transformation groups. A trivial example of this principle is the
invariance of the theory with respect to Euclidean motions. A less
trivial but well-known example is the discrete analog for the
local Gaussian curvature defined as the angle defect $ G(v)=2\pi
-\sum \alpha_i, $ at a vertex $v$ of a polyhedral surface. Here
the $\alpha_i$ are the angles of all polygonal faces (see
Fig.~\ref{f.dEnergy}) of the surface at vertex $v$. The discrete
Gaussian curvature $G(v)$ defined in this way is preserved under
isometries, which is a discrete version of the theoremum egregium
of Gauss.

\begin{figure}[h]
\begin{center}
\hspace{1cm}\parbox[c]{0.25\textwidth}{{\hspace{0.3pt}\includegraphics[width=0.25\textwidth]{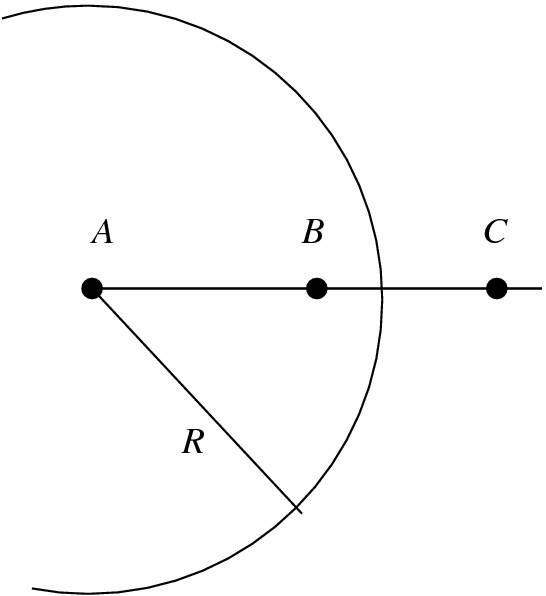}}}
\hfill
\parbox[c]{0.65\textwidth}{{\hspace{0.5pt}\includegraphics[width=0.65\textwidth]{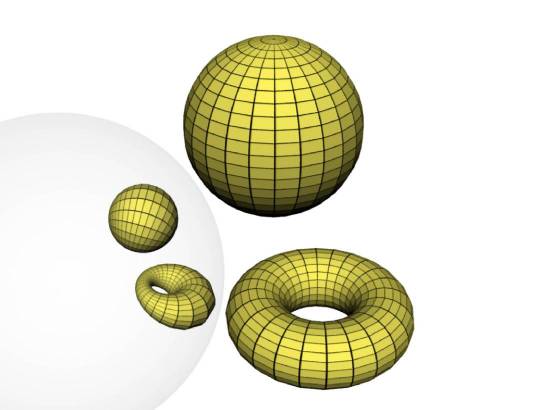}}}
\end{center}
\caption{Inversion $B \mapsto C$ in a sphere, $\mid AB\mid \mid
AC\mid =R^2$. A sphere and a torus of revolution and their
inversions in a sphere: spheres are mapped to spheres.
\label{f.inversion}}
\end{figure}

In these lectures, we focus on surface geometries invariant under
M\"obius transformations. Recall that M\"obius transformations
form a finite-dimensional Lie group generated by inversions in
spheres (see Fig.~\ref{f.inversion}).
M\"obius transformations can be also thought as compositions of
translations, rotations, homotheties and inversions in spheres.
Alternatively, in dimensions $n\ge 3$ M\"obius transformations can
be characterized as conformal transformations: Due to Liouville's
theorem any conformal mapping $F:U\to V$ between two open subsets
$U,V\subset {\mathbb R}^n, n\ge 3$ is a M\"obius transformation.

Many important geometric notions and properties are known to be
preserved by M\"obius transformations. The list includes in
particular:
\begin{itemize}
    \item spheres of any dimension, in particular circles (planes and straight
    lines are treated as infinite spheres and circles),
    \item intersection angles between spheres (and circles),
    \item curvature line parametrization,
    \item conformal parametrization,
    \item isothermic parametrization ( conformal curvature line parametrization),
    \item the Willmore functional (see Section~\ref{s.willmore}).
\end{itemize}

For discretization of M\"obius-invariant notions it is natural to
use M\"obius-invariant building blocks. This observation leads us
to the conclusion that the discrete conformal or curvature line
parametrizations of surfaces and the discrete Willmore functional
should be formulated in terms of circles and spheres.

\section{Discrete Willmore Energy}
\label{s.willmore}

The Willmore functional \cite{Willmore} for a smooth surface $\mathcal S$ in
3-dimensional Euclidean space is
$$
{\mathcal W(S)}=\frac{1}{4}\int_{\mathcal S} (k_1-k_2)^2 dA=\int_{\mathcal S} H^2
dA- \int_{\mathcal S} K dA.
$$
Here $dA$ is the area element, $k_1$ and $k_2$ the principal
curvatures, $H=\frac{1}{2}(k_1+k_2)$ the mean curvature, and
$K=k_1 k_2$ the Gaussian curvature of the surface.

Let us mention two important properties of the Willmore energy:
\begin{itemize}
    \item  ${\mathcal W(S)}\ge 0$ and ${\mathcal W(S)}=0$ if and only if
    $\mathcal S$ is a round sphere.
    \item  $\mathcal W(S)$ (and the integrand $(k_1-k_2)^2 dA$) is M\"obius
    invariant \cite{Blaschke,Willmore}
\end{itemize}
Whereas the first claim almost immediately follows from the
definition, the second is a non-trivial property. Observe that for
closed surfaces $\mathcal W(S)$ and $\int_{\mathcal S} H^2 dA$
differ by a topological invariant $\int K dA=2\pi \chi(S)$. We
prefer the definition of $\mathcal W(S)$ with a M\"obius invariant
integrand.

Observe that minimization of the Willmore energy $\mathcal W$
seeks to make the surface ``as round as possible''. This property
and the M\"obius invariance are two principal points of the
geometric discretization of the Willmore energy suggested in
\cite{BobenkoWillmore}. In this section we present the main
results of \cite{BobenkoWillmore} with complete derivations, some
of which were omitted there.

\subsection{Discrete Willmore functional for simplicial surfaces}

Let $S$ be a simplicial surface in 3-dimensional Euclidean space
with vertex set $V$, edges $E$ and (triangular) faces $F$. We
define the discrete Willmore energy of $S$ using the circumcircles
of its faces. Each (internal) edge $e\in E$ is incident to two
triangles. A consistent orientation of the triangles naturally
induces an orientation of the corresponding circumcircles. Let
$\beta(e)$ be the external intersection angle of the circumcircles
of the triangles sharing $e$, meaning the angle between the
tangent vectors of the oriented circumcircles (at either
intersection point).

\begin{definition}
The local discrete Willmore energy at a vertex $v$ is the sum
$$
W(v)=\sum_{e\ni v} \beta(e)-2\pi.
$$
over all edges incident to $v$. The discrete Willmore energy of a compact
simplicial surface $S$ without boundary is the sum over all vertices
$$
W(S)=\dfrac{1}{2}\sum_{v\in V}W(v)=\sum_{e\in E} \beta(e)-\pi\mid V \mid.
$$
Here $|V|$ is the number of vertices of $S$.
\end{definition}

\begin{figure}[h]
\begin{center}
\parbox[c]{0.4\textwidth}{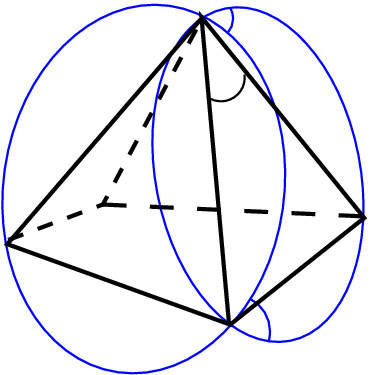}
\parbox[c]{0.3\textwidth}{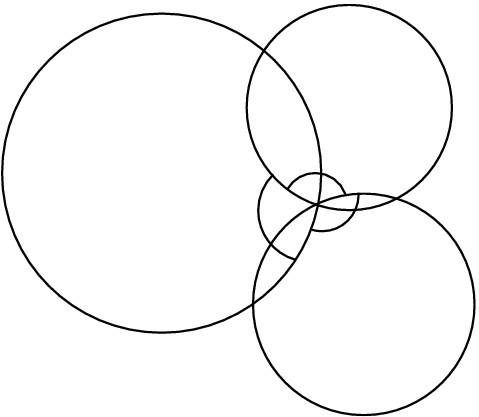}
\end{center}
\caption{Definition of discrete Willmore energy
\label{f.dEnergy}}
\end{figure}

Figure~\ref{f.dEnergy} presents two neighboring circles with their external
intersection angle $\beta_i$ as well as a view ``from the top'' at a vertex $v$
showing all $n$ circumcircles passing through $v$ with the corresponding
intersection angles $\beta_1,\ldots,\beta_n$. For simplicity we will consider only
simplicial surfaces without boundary.

The energy $W(S)$ is obviously invariant with respect to M\"obius transformations.

The star $S(v)$ of the vertex $v$ is the subcomplex of $S$
consisting of the triangles incident with $v$. The vertices of
$S(v)$ are $v$ and all its neighbors. We call $S(v)$ convex if for
each of its faces $f\in F(S(v))$ the star $S(v)$ lies to one side
of the plane of $f$ and strictly convex if the intersection of
$S(v)$ with the plane of $f$ is $f$ itself.

\begin{proposition} \label{p.non-negative}
The conformal energy is nonnegative
$$W(v)\ge 0,$$
and vanishes if and only if the star $S(v)$ is convex and all its vertices lie on a
common sphere.
\end{proposition}

\begin{figure}[h]
\begin{center}
\parbox[c]{0.4\textwidth}{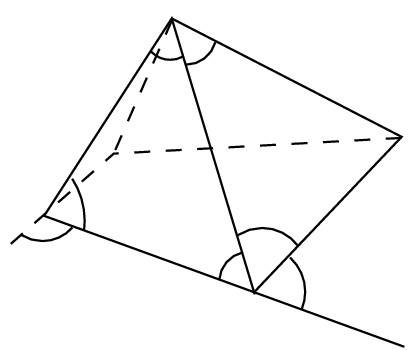}
\end{center}
\caption{Proof of Lemma \ref{l.alpha-beta} \label{f.alpha-beta}}
\end{figure}

The proof of this proposition is based on an elementary lemma.
\begin{lemma} \label{l.alpha-beta}
Let $\mathcal P$ be a (not necessarily planar) $n$-gon with external angles
$\beta_i$. Choose a point $P$ and connect it to all vertices of $\mathcal P$. Let
$\alpha_i$ be the angles of the triangles at the tip $P$ of the pyramid thus
obtained (see Figure~\ref{f.alpha-beta}). Then
$$
\sum_{i=1}^n \beta_i\ge \sum_{i=1}^n \alpha_i,
$$
and equality holds if and only if $\mathcal P$ is planar and convex
and the vertex $P$ lies inside $\mathcal P$.
\end{lemma}

\begin{proof}
Denote by $\gamma_i$ and $\delta_i$ the angles of the triangles at
the vertices of $\mathcal P$, as in Figure~\ref{f.alpha-beta}. The
claim of Lemma \ref{l.alpha-beta} follows from adding over all
$i=1,\ldots,n$  the two obvious relations
\begin{eqnarray*}
\beta_{i+1}&\ge& \pi -(\gamma_{i+1}+\delta_i)\\
\alpha_i&=&\pi -(\gamma_{i}+\delta_i).
\end{eqnarray*}
All inequalities become equalities only in the case when $\mathcal P$ is planar,
convex and contains $P$.
\end{proof}

For $P$ in the convex hull of $\mathcal P$ we have $\sum \alpha_i
\ge 2\pi$. As a corollary we obtain a polygonal version of
Fenchel's theorem \cite{Fenchel}:

\begin{corollary}  \label{c.non-negative}
$$
\sum_{i=1}^n \beta_i\ge 2\pi.
$$
\end{corollary}

\noindent{\it Proof of Proposition \ref{p.non-negative}}. The
claim of Proposition \ref{p.non-negative} is invariant with
respect to M\"obius transformations. Applying a M\"obius
transformation $M$ that maps the vertex $v$ to infinity,
$M(v)=\infty$, we make all circles passing through $v$ into
straight lines and arrive at the geometry shown in
Figure~\ref{f.alpha-beta} with $P=M(\infty)$. Now the result
follows immediately from Corollary \ref{c.non-negative}.

\begin{theorem}
Let $S$ be a compact simplicial surface without boundary. Then
$$W(S)\ge 0,$$
and equality holds if and only if $S$ is a convex polyhedron
inscribed in a sphere, i.e. a Delaunay triangulation of a sphere.
\end{theorem}
\begin{proof}
 Only the second statement needs to be proven. By Proposition
\ref{p.non-negative},  the equality $W(S)=0$ implies that the star
of each vertex of $S$ is convex (but not necessarily strictly
convex). Deleting the edges that separate triangles lying in a
common plane, one obtains a polyhedral surface $S_P$ with circular
faces and all strictly convex vertices and edges. Proposition
\ref{p.non-negative} implies that for every vertex $v$ there
exists a sphere $S_v$ with all vertices of the star $S(v)$ lying
on it. For any edge $(v_1,v_2)$ of $S_P$ two neighboring spheres
$S_{v_1}$ and $S_{v_2}$ share two different circles of their
common faces. This implies $S_{v_1}=S_{v_2}$ and finally the
coincidence of all the spheres $S_v$.
\end{proof}

\subsection{Noninscribable polyhedra}

The minimization of the conformal energy for simplicial spheres is
related to a classical result of Steinitz \cite{Steinitz}, who
showed that there exist abstract simplicial 3-polytopes without
geometric realizations as convex polytopes with all vertices on a
common sphere. We call these combinatorial types noninscribable.

Let $S$ be a simplicial sphere with vertices colored in black and white. Denote the
sets of white and black vertices by $V_w$ and $V_b$ respectively, $V=V_w\cup V_b$.
Assume that there are no edges connecting two white vertices and denote the sets of
the edges connecting white and black vertices and two black vertices by $E_{wb}$
and $E_{bb}$ respectively, $E=E_{wb}\cup E_{bb}$. The sum of the local discrete
Willmore energies over all white vertices can be represented as
$$
\sum_{v\in V_w} W(v)=\sum_{e\in E_{wb}}\beta(e)-2\pi |V_w|.
$$
Its nonnegativity yields $\sum_{e\in E_{wb}}\beta(e)\ge 2\pi
|V_w|$. For the discrete Willmore energy of $S$ this implies
\begin{equation}    \label{e.steinitz}
W(S)=\sum_{e\in E_{wb}}\beta(e)+\sum_{e\in E_{bb}}\beta(e)-
\pi(|V_w| +|V_b|)\ge \pi(|V_w| -|V_b|).
\end{equation}
The equality here holds if and only if $\beta(e)=0$ for all $e\in
E_{bb}$ and the star of any white vertices is convex, with
vertices lying on a common sphere. We come to the conclusion that
the polyhedra of this combinatorial type with $|V_w| > |V_b|$ have
positive Willmore energy and thus cannot be realized as convex
polyhedra all of whose vertices belong to a sphere. These are
exactly the noninscribable examples of Steinitz \cite{Gruenbaum}.

\begin{figure}[h]
\begin{center}
\includegraphics[width=0.25\textwidth]{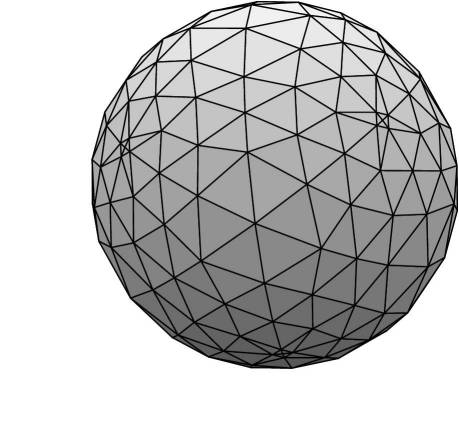}\hspace{3cm}
\includegraphics[width=0.25\textwidth]{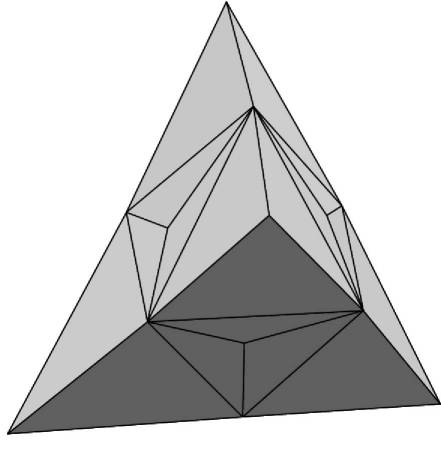}
\end{center}
\caption{Discrete Willmore spheres of inscribable ($W=0$) and non-inscribable
($W>0$) types. \label{f.minimizing}}
\end{figure}

One such example is presented in Figure~\ref{f.minimizing}. Here
the centers of the edges of the tetrahedron are black and all
other vertices are white, so $|V_w|=8, |V_b|=6$. The estimate
(\ref{e.steinitz}) implies that the discrete Willmore energy of
any polyhedron of this type is at least $2\pi$. The polyhedra with
energy equal to $2\pi$ are constructed as follows. Take a
tetrahedron, color its vertices white and chose one black vertex
per edge. Draw circles through each white vertex and its two black
neighbors. We get three circles on each face. Due to Miquel's
theorem (see Fig.~\ref{f.Miquel}) these three circles intersect at
one point. Color this new vertex white. Connect it by edges to all
black vertices of the triangle and connect pairwise the black
vertices of the original faces of the tetrahedron. The constructed
polyhedron has $W=2\pi$.

To construct further polyhedra with $|V_w|>|V_b|$, take a
polyhedron $\hat{P}$ whose number of faces is greater then the
number of vertices $|\hat{F}|> |\hat{V}|$. Color all the vertices
black, add white vertices at the faces and connect them to all
black vertices of a face. This yields a polyhedron with
$|V_w|=|\hat{F}| > |V_b|=|\hat{V}|$. Hodgson, Rivin and Smith
\cite{HodgsonRS} have found a characterization of inscribable
combinatorial types, based on a transfer to the Klein model of
hyperbolic 3-space. Their method is related to the methods of
construction of discrete minimal surfaces in
Section~\ref{s.minimal}.

The example in Figure~\ref{f.minimizing}~(right) is one of the few
for which the minimum of the discrete Willmore energy can be found
by elementary methods. Generally this is a very appealing (but
probably difficult) problem of discrete differential geometry (see
the discussion in \cite{BobenkoWillmore}).

Complete understanding of noninscribable simplicial spheres is an
interesting mathematical problem. However the phenomenon of
existence of such spheres might be seen as a problem in using of
the discrete Willmore functional for applications in computer
graphics, such as the fairing of surfaces. Fortunately the problem
disappears after just one refinement step: all simplicial spheres
become inscribable. Let $\bf S$ be an abstract simplicial sphere.
Define its refinement $\bf S_R$ as follows: split every edge of
$\bf S$ in two by inserting additional vertices, and connect these
new vertices sharing a face of $\bf S$ by additional edges ($1\to
4$ refinement, see Figure~\ref{f.smooth}~(left)).

\begin{proposition}
\label{p.refine} The refined simplicial sphere $\bf S_R$ is inscribable, and thus
there exists a polyhedron $S_R$ with the combinatorics of $\bf S_R$ and $W(S_R)=0$.
\end{proposition}

\begin{proof}
Koebe's theorem (see Theorem~\ref{t.polyKoebe},
Section~\ref{s.minimal}) states that every abstract simplicial
sphere $\bf S$ can be realized as a convex polyhedron $S$ all of
whose edges touch a common sphere $S^2$. Starting with this
realization $S$ it is easy to construct a geometric realization
$S_R$ of the refinement $\bf S_R$ inscribed in $S^2$. Indeed,
choose the touching points of the edges of $S$ with $S^2$ as the
additional vertices of $S_R$ and project the original vertices of
$S$ (which lie outside of the sphere $S^2$) to $S^2$. One obtains
a convex simplicial polyhedron $S_R$ inscribed in $S^2$.
\end{proof}

\subsection{Computation of the energy}

For derivation of some formulas it will be convenient to use the
language of quaternions. Let $\{{\bf 1},{\bf i},{\bf j},{\bf k}\}$
be the standard basis
$$
{\bf ij}={\bf k},\quad {\bf ij}={\bf k},\quad {\bf ij}={\bf
k},\quad {\bf ii}={\bf jj}={\bf kk}=-{\bf 1}
$$
of the quaternion algebra ${\mathbb H}$. A quaternion $q=q_0{\bf
1}+q_1{\bf i}+q_2{\bf j}+q_3{\bf k}$ is decomposed in its real
part ${\rm Re}\ q:=q_0\in{\mathbb R}$ and imaginary part ${\rm
Im}\ q:=q_1{\bf i}+q_2{\bf j}+q_3{\bf k}\in{\rm Im}\ {\mathbb H}$.
The absolute value of $q$ is $|q| := q_0^2+q_1^2+q_2^2+q_3^2$.

We identify vectors in ${\mathbb R}^3$ with imaginary quaternions
$$
v=(v_1,v_2, v_3)\in{\mathbb R}^3\quad \longleftrightarrow \quad
v=v_1{\bf i}+v_2{\bf j}+v_3{\bf k}\in{\rm Im}\ {\mathbb H}
$$
and do not distinguish them in our notation. For the quaternionic
product this implies
\begin{equation}   \label{product}
vw=-\langle v,w\rangle+v\times w,
\end{equation}
where $\langle v,w\rangle$ and $v\times w$ are the scalar and vector products in
${\mathbb R}^3$.

\begin{definition}
Let $x_1,x_2,x_3,x_4\in {\mathbb R}^3\cong{\rm Im}\ {\mathbb H}$ be points in
3-dimensional Euclidean space. The quaternion
$$
q(x_1,x_2,x_3,x_4):=(x_1-x_2)(x_2-x_3)^{-1}(x_3-x_4)(x_4-x_1)^{-1}
$$
is called the cross-ratio of $x_1,x_2,x_3,x_4$.
\end{definition}

The cross-ratio is quite useful due to its M\"obius properties:
\begin{lemma}
The absolute value and real part of the cross-ratio
$q(x_1,x_2,x_3,x_4)$ (or equivalently $|q|$ and $| {\rm Im}\ q|$)
are preserved by M\"obius transformations. The quadrilateral
$x_1,x_2,x_3,x_4$ is circular if and only if
$q(x_1,x_2,x_3,x_4)\in{\mathbb R}$.
\end{lemma}

\begin{figure}[h]
\begin{center}
\parbox[c]{0.4\textwidth}{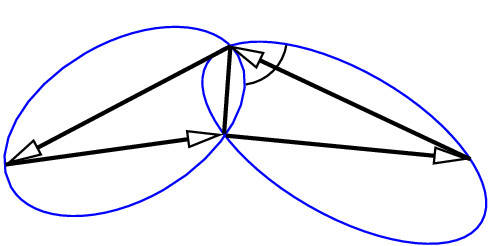}
\end{center}
\caption{Formula for the angle between circumcircles \label{f.formula}}
\end{figure}

Consider two triangles with a common edge. Let $a, b, c, d\in {{\mathbb R}^3}$ be
their other edges, oriented as in Fig.\ref{f.formula}.

\begin{proposition}
The external angle $\beta\in [0,\pi ]$ between the circumcircles
of the triangles in Figure~\ref{f.formula} is given by any of the
equivalent formulas:
\begin{eqnarray}
\label{cos-beta}
\cos(\beta)&=&-\dfrac{{\rm Re}\ q}{\mid q\mid}=
            -\dfrac{{\rm Re}\ (abcd)}{\mid abcd\mid}=  \nonumber\\
           &=&\dfrac{\langle a,c\rangle\langle b,d\rangle -
           \langle a,b\rangle\langle c,d\rangle -\langle b,c\rangle\langle d,a\rangle }
           {\mid a\mid \mid b\mid \mid c\mid \mid d\mid}.
\end{eqnarray}
Here $q=ab^{-1}cd^{-1}$ is the cross-ratio of the quadrilateral.
\end{proposition}
\begin{proof} Since ${\rm Re}\ q$, $|q|$ and $\beta$ are
M\"obius-invariant, it is enough to prove the first formula for
the planar case $a,b,c,d\in {\mathbb C}$, mapping all four
vertices to a plane by a M\"obius transformation. In this case $q$
becomes the classical complex cross-ratio. Considering the
arguments $a,b,c,d\in {\mathbb C}$ one easily arrives at
$\beta=\pi-\arg q$. The second representation follows from the
identity $b^{-1}=-b/|b|$ for imaginary quaternions. Finally
applying (\ref{product}) we obtain
\begin{eqnarray*}
{\rm Re}\ (abcd)=\langle a,b\rangle \langle c,d\rangle -
\langle a\times b,c\times d\rangle =\\
\langle a,b\rangle\langle c,d\rangle +\langle b,c\rangle\langle d,a\rangle -
\langle a,c\rangle\langle b,d\rangle .
\end{eqnarray*}
\end{proof}

\subsection{Smooth limit}

\begin{figure}[h]
\begin{center}
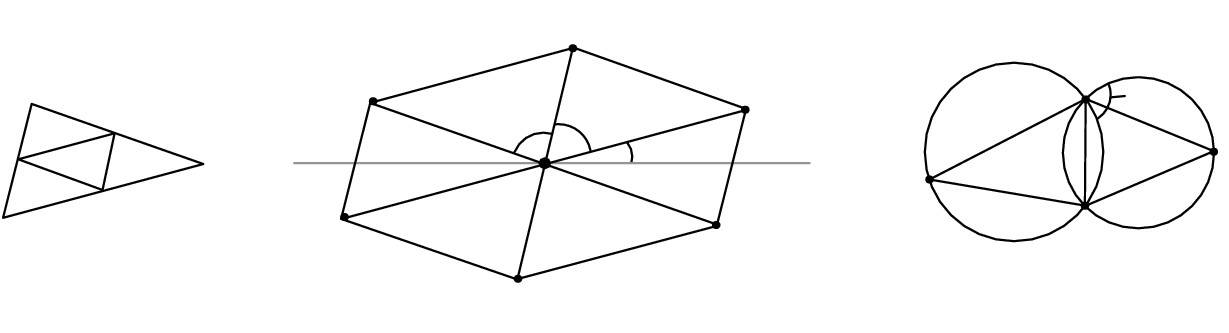
\end{center}
\caption{Smooth limit of the discrete Willmore energy. {\em Left:} The $1\to 4 $
refinement. {\em Middle:} An infinitesimal hexagon in the parameter plane with a
(horizontal) curvature line. {\em Right:} The $\beta$-angle corresponding to two
neighboring triangles in ${\mathbb R}^3$.
 \label{f.smooth}}
\end{figure}

The discrete energy $W$ is not only a discrete analogue of the
Willmore energy. In this section we show that it approximates the
smooth Willmore energy, although the smooth limit is very
sensitive to the refinement method and should be chosen in a
special way. We consider a special infinitesimal triangulation
which can be obtained in the limit of $1\to 4$ refinements (see
Fig.~\ref{f.smooth} (left)) of a triangulation of a smooth
surface. Intuitively it is clear that in the limit one has a
regular triangulation such that almost every vertex is of valence
6 and neighboring triangles are congruent up to sufficiently high
order of $\epsilon$ ($\epsilon$ is the order of the distances
between neighboring vertices).

We start with a comparison of the discrete and smooth Willmore
energies for an important modelling example. Consider a
neighborhood of a vertex $v\in {\mathcal S}$, and represent the
smooth surface locally as a graph over the tangent plane at $v$:
$$
{\mathbb R}^2\ni(x,y)\mapsto f(x,y)=\Big(x, y, {1\over 2}(k_1
x^2+k_2 y^2) + o(x^2+y^2)\Big)\in{\mathbb R}^3,\quad (x,y)\to 0.
$$
Here $x,y$ are the curvature directions and $k_1,k_2$ are the
principal curvatures at $v$. Let the vertices $(0,0), a=(a_1,a_2)$
and $b=(b_1,b_2)$ in the parameter plane form an acute triangle.
Consider the infinitesimal hexagon with vertices $\epsilon a,
\epsilon b, \epsilon c, -\epsilon a, -\epsilon b, -\epsilon c$,
(see Figure~\ref{f.smooth} (middle)), with $b=a+c$. The
coordinates of the corresponding points on the smooth surface are
\begin{eqnarray*}
f(\pm\epsilon a)&=&\epsilon(\pm a_1,\pm a_2,\epsilon r_a+o(\epsilon)),\\
f(\pm\epsilon c)&=&\epsilon(\pm c_1,\pm c_2,\epsilon r_c+o(\epsilon)),\\
f(\pm\epsilon b)&=& (f(\pm\epsilon a)+ f(\pm\epsilon c))+\epsilon^2 R,\quad R=(0,0,
r+o(\epsilon)),
\end{eqnarray*}
where
$$
r_a={1\over 2}(k_1 a_1^2+k_2 a_2^2),\quad r_c={1\over 2}(k_1 c_1^2+k_2 c_2^2),
\quad r=(k_1 a_1 c_1+k_2 a_2 c_2)
$$
and $a=(a_1,a_2), c=(c_1,c_2)$.

We will compare the discrete Willmore energy $W$ of the simplicial
surface comprised by the vertices $f(\epsilon a),\ldots,
f(-\epsilon c)$ of the hexagonal star with the classical Willmore
energy ${\mathcal W}$ of the corresponding part of the smooth
surface $\mathcal S$. Some computations are required for this.
Denote by $\epsilon A=f(\epsilon a),\epsilon B=f(\epsilon
b),\epsilon C=f(\epsilon c)$ the vertices of two corresponding
triangles (as in Figure~\ref{f.smooth} (right)), and also by $|a|$
the length of $a$ and by $\langle a,c\rangle=a_1c_1+a_2c_2$ the
corresponding scalar product.
\begin{lemma}
The external angle $\beta(\epsilon)$ between the circumcircles of
the triangles with the vertices $(0,A,B)$ and $(0,B,C)$ (as in
Figure~\ref{f.smooth} (right)) is given by
\begin{eqnarray}        \label{beta-epsilon}
\beta(\epsilon)=\beta(0)+ w(b)+o(\epsilon^2), \quad \epsilon\to 0,\quad
w(b)=\epsilon^2\dfrac{g\cos\beta(0)-h}{|a|^2|c|^2\sin\beta(0)}.
\end{eqnarray}
Here $\beta(0)$ is the external angle of the circumcircles of the triangles
$(0,a,b)$ and $(0,b,c)$ in the plane, and
\begin{eqnarray*}
g&=&|a|^2r_c(r+r_c)+|c|^2r_a(r+r_a)+{r^2\over 2}(|a|^2+|c|^2),\\
h&=&|a|^2r_c(r+r_c)+|c|^2r_a(r+r_a)- \langle a,c\rangle (r+2r_a)(r+2r_c).
\end{eqnarray*}
\end{lemma}
\begin{proof}
Formula (\ref{cos-beta}) with $a=-C,b=A,c=C+\epsilon R, d=-A-\epsilon R$ yields
$$
\cos\beta=\dfrac{\langle C,C+\epsilon R\rangle\langle A,A+\epsilon R\rangle-
            \langle A,C\rangle\langle A+\epsilon R,C+\epsilon R\rangle-
            \langle A,C+\epsilon R\rangle
            \langle A+\epsilon R,C\rangle}{|A||C||A+\epsilon R||C+\epsilon
            R|},
$$
where $|A|$ is the length of $A$. Substituting the expressions for
$A,C,R$ we see that the term of order $\epsilon$ of the numerator
vanishes, and we obtain for the numerator
$$
|a|^2|c|^2-2\langle a,c\rangle^2+\epsilon^2 h +o(\epsilon^2).
$$
For the terms in the denominator we get
$$
|A|=|a|\left(1+{r_a^2\over
2|a|^2}\epsilon^2+o(\epsilon^2)\right),\quad |A+\epsilon
R|=|a|\left(1+{(r+r_a)^2\over
2|a|^2}\epsilon^2+o(\epsilon^2)\right)
$$
and similar expression for $|C|$ and $|C+\epsilon R|$. Substituting this to the
formula for $\cos\beta$ we obtain
$$
\cos\beta=1-2\left(\dfrac{\langle a,c\rangle}{|a||c|}\right)^2+
\dfrac{\epsilon^2}{|a|^2|c|^2}\left(h-g(1-2\left(\dfrac{\langle
a,c\rangle}{|a||c|}\right)^2)   \right)+o(\epsilon^2).
$$
Observe that this formula can be read as
$$
\cos\beta(\epsilon)=\cos\beta(0)+\dfrac{\epsilon^2}{|a|^2|c|^2}\left(h-g\cos\beta(0)
\right)+o(\epsilon^2),
$$
which implies the asymptotic (\ref{beta-epsilon}).
\end{proof}
The term $w(b)$ is in fact the part of the discrete Willmore energy of the vertex
$v$ coming from the edge $b$. Indeed the sum of the angles $\beta(0)$ over all 6
edges meeting at $v$ is $2\pi$. Denote by $w(a)$ and $w(c)$ the parts of the
discrete Willmore energy corresponding to the edges $a$ and $c$. Observe that for
the opposite edges (for example $a$ and $-a$) the terms $w$ coincide. Denote
$W_\epsilon(v)$ the discrete Willmore energy of the simplicial hexagon we consider.
We have
$$
W_\epsilon(v)=(w(a)+w(b)+w(c))+o(\epsilon^2).
$$
On the other hand the part of the classical Willmore functional corresponding to
the vertex $v$ is
$$
{\mathcal W}_\epsilon(v)={1\over 4}(k_1-k_2)^2 S +0(\epsilon^2),
$$
where the  area $S$ is one third of the area of the hexagon or, equivalently, twice
the area of one of the triangles in the parameter domain
$$
S=\epsilon^2|a||c|\sin\gamma.
$$
Here $\gamma$ is the angle between the vectors $a$ and $c$. An elementary geometric
consideration implies
\begin{equation}            \label{beta-gamma}
\beta(0)=2\gamma-\pi.
\end{equation}
We are interested in the quotient $W_\epsilon/{\mathcal W}_\epsilon$ which is
obviously scale invariant. Let us normalize $|a|=1$ and parametrize the triangles
by the angles between the edges and by the angle to the curvature line (see
Fig.~\ref{f.smooth} (middle)).
\begin{eqnarray}
(a_1,a_2)&=&(\cos\phi_1,\sin\phi_1),            \label{ac-phi}\\
(c_1,c_2)&=&(\dfrac{\sin\phi_2}{\sin\phi_3}\cos(\phi_1+\phi_2+\phi_3),
            \dfrac{\sin\phi_2}{\sin\phi_3}\sin(\phi_1+\phi_2+\phi_3)).\nonumber
\end{eqnarray}
The moduli space of the regular lattices of acute triangles is described as follows
$$
\Phi= \{ \phi=(\phi_1,\phi_2,\phi_3)\in {\mathbb R}^3| 0\le\phi_1<{\pi\over 2},\
0<\phi_2<{\pi\over 2},\ 0<\phi_3<{\pi\over 2},\ {\pi\over 2}<\phi_2+\phi_3\}.
$$
\begin{proposition}
\label{p.limit}
The limit of the quotient of the discrete and
smooth Willmore energies
$$
Q(\phi):=\lim_{\epsilon\to 0}\dfrac{W_\epsilon(v)}{{\mathcal W}_\epsilon(v)}
$$
is independent of the curvatures of the surface and depends on the
geometry of the triangulation only. It is
\begin{equation}                \label{Q}
Q(\phi)=1-\dfrac{(\cos2\phi_1\cos\phi_3+\cos(2\phi_1+2\phi_2+\phi_3))^2+
(\sin2\phi_1\cos\phi_3)^2}{4 \cos\phi_2\cos\phi_3\cos(\phi_2+\phi_3)},
\end{equation}
and we have $Q>1$. The infimum $\inf_\Phi Q(\phi)=1$ corresponds
to one the cases when two of the three lattice vectors $a,b,c$ are
in the principal curvature directions
\begin{itemize}
    \item $\phi_1=0, \phi_2+\phi_3\to{\pi\over 2}$,
    \item $\phi_1=0, \phi_2\to{\pi\over 2}$,
    \item $\phi_1+\phi_2={\pi\over 2},\phi_3\to{\pi\over 2}$.
\end{itemize}
\end{proposition}
\begin{proof} is based on a direct but rather involved computation. We have
used the {\em Mathematica} computer algebra system for some of the
computations. Introduce $$ \tilde{w}:= {4w\over (k_1-k_2)^2S}.
$$
This gives in particular
$$
\tilde{w}(b)=2\dfrac{h+g(2\cos^2\gamma-1)}{(k_1-k_2)^2|a|^3|c|^3\cos\gamma\sin^2\gamma}=
2\dfrac{h+g(2{\langle a,c\rangle^2\over |a|^2|c|^2}-1)}{(k_1-k_2)^2\langle
a,c\rangle(|a|^2|c|^2-\langle a,c\rangle^2)}.
$$
Here we have used the relation (\ref{beta-gamma}) between $\beta(0)$ and $\gamma$.
In the sum over the edges $Q=\tilde{w}(a)+\tilde{w}(b)+\tilde{w}(c)$ the curvatures
$k_1,k_2$ disappear and we get $Q$ in terms of the coordinates of $a$ and $c$:
\begin{eqnarray*}
Q&=&2((a_1^2c_2^2+a_2^2c_1^2)(a_1c_1+a_2c_2)+ a_1^2c_1^2(a_2^2+c_2^2)+
a_2^2c_2^2(a_1^2+c_1^2)+
\\& &2a_1a_2c_1c_2((a_1+c_1)^2+(a_2+c_2)^2) )/ \\
& &\left((a_1c_1+a_2c_2)
(a_1(a_1+c_1)+a_2(a_2+c_2))((a_1+c_1)c_1+(a_2+c_2)c_2)\right).
\end{eqnarray*}
Substituting the angle representation (\ref{ac-phi}) we obtain
$$
Q= \dfrac{\sin2\phi_1\sin2(\phi_1+\phi_2)+
2\cos\phi_2\sin(2\phi_1+\phi_2)\sin2(\phi_1+\phi_2+\phi_3)}{4
\cos\phi_2\cos\phi_3\cos(\phi_2+\phi_3)}.
$$
One can check that this formula is equivalent to (\ref{Q}). Since the denominator
in (\ref{Q}) on the space $\Phi$ is always negative we have $Q>1$. The identity
$Q=1$ holds only if both terms in the nominator of (\ref{Q}) vanish. This leads
exactly to the cases indicated in the proposition when the lattice vectors are
directed along the curvature lines. Indeed the vanishing of the second term in the
nominator implies either $\phi_1=0$ or $\phi_3\to {\pi\over 2}$. Vanishing of the
first term in the nominator with $\phi_1=0$ implies $\phi_2\to {\pi\over 2}$ or
$\phi_2+\phi_3\to {\pi\over 2}$. Similarly in the limit $\phi_3\to {\pi\over 2}$
the vanishing of
$$
{(\cos2\phi_1\cos\phi_3+\cos(2\phi_1+2\phi_2+\phi_3))^2}/{\cos\phi_3}
$$
implies $\phi_1+\phi_2={\pi\over 2}$. One can check that in all these cases
$Q(\phi)\to 1$.
\end{proof}
Note that for the infinitesimal equilateral triangular lattice
$\phi_2=\phi_3={\pi\over 3}$ the result is independent of the orientation $\phi_1$
with respect to the curvature directions, and the discrete Willmore energy is in
the limit Q=3/2 times larger than the smooth one.

Finally, we come to the following conclusion.
\begin{theorem}             \label{t.smooth}
Let $\mathcal S$ be a smooth surface with Willmore energy
${\mathcal W}({\mathcal S})$. Consider a simplicial surface
$S_\epsilon$ such that its vertices lie on $\mathcal S$ and are of
degree 6, the distances between the neighboring vertices are of
order $\epsilon$, and the neighboring triangles of $S_\epsilon$
meeting at a vertex are congruent up to order $\epsilon^3$ (i.e.
the lengths of the corresponding edges differ by terms of order at
most $\epsilon^4$), and they build elementary hexagons the lengths
of whose opposite edges differ by terms of order at most
$\epsilon^4$. Then the limit of the discrete Willmore energy is
bounded from below by the classical Willmore energy
\begin{equation}                    \label{W>W}
\lim_{\epsilon\to 0}W(S_\epsilon)\ge {\mathcal W}({\mathcal S}).
\end{equation}
Moreover equality in (\ref{W>W}) holds if $S_\epsilon$ is a
regular triangulation of an infinitesimal curvature line net of
$\mathcal S$, i.e. the vertices of $S_\epsilon$ are at the
vertices of a curvature line net of $\mathcal S$.
\end{theorem}

\begin{proof}
Consider an elementary hexagon of $S_\epsilon$. Its projection to
the tangent plane of the central vertex is a hexagon which can be
obtained from the modelling one considered in
Proposition~\ref{p.limit} by a perturbation of vertices of order
$o(\epsilon^3)$. Such perturbations contribute to the terms of
order $o(\epsilon^2)$ of the discrete Willmore energy. The latter
are irrelevant for the considerations of
Proposition~\ref{p.limit}.
\end{proof}

Possibly minimization of the discrete Willmore energy with the
vertices constrained to lie on $S$ could be used for computation
of a curvature-line net.

\subsection{Bending energy for simplicial surfaces}

An accurate model for bending of discrete surfaces is important
for modeling in computer graphics.

The bending energy of smooth thin shells is given by the integral \cite{Grinspun}
$$
E=\int (H-H_0)^2 dA,
$$
where $H_0$ and $H$ are the mean curvatures of the original and deformed surface
respectively. For $H_0=0$ it reduces to the Willmore energy.

\begin{figure}[h]
\begin{center}
\parbox[c]{0.4\textwidth}
{\begin{picture}(0,0)%
\includegraphics{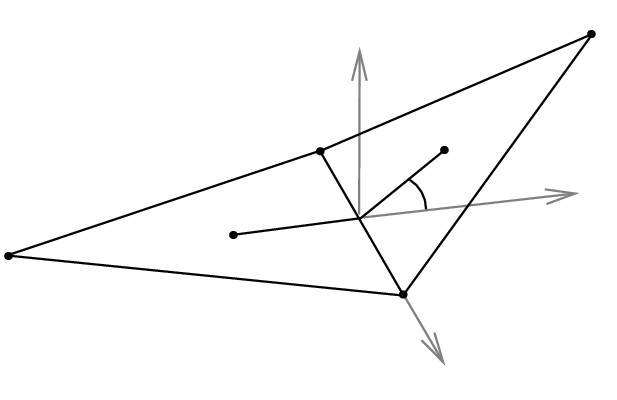}%
\end{picture}%
\setlength{\unitlength}{4144sp}%
\begingroup\makeatletter\ifx\SetFigFont\undefined%
\gdef\SetFigFont#1#2#3#4#5{%
  \reset@font\fontsize{#1}{#2pt}%
  \fontfamily{#3}\fontseries{#4}\fontshape{#5}%
  \selectfont}%
\fi\endgroup%
\begin{picture}(2341,1515)(348,-1702)
\put(1984,-942){\makebox(0,0)[lb]{\smash{{\SetFigFont{9}{10.8}{\familydefault}{\color[rgb]{0,0,0}$\theta$}%
}}}}
\put(2598,-952){\makebox(0,0)[lb]{\smash{{\SetFigFont{9}{10.8}{\familydefault}{\color[rgb]{0,0,0}2}%
}}}}
\put(1442,-683){\makebox(0,0)[lb]{\smash{{\SetFigFont{9}{10.8}{\familydefault}{\color[rgb]{0,0,0}$X_4$}%
}}}}
\put(1024,-1109){\makebox(0,0)[lb]{\smash{{\SetFigFont{9}{10.8}{\familydefault}{\color[rgb]{0,0,0}$X_1$}%
}}}}
\put(1687,-307){\makebox(0,0)[lb]{\smash{{\SetFigFont{9}{10.8}{\familydefault}{\color[rgb]{0,0,0}3}%
}}}}
\put(2058,-1687){\makebox(0,0)[lb]{\smash{{\SetFigFont{9}{10.8}{\familydefault}{\color[rgb]{0,0,0}1}%
}}}}
\put(1952,-1337){\makebox(0,0)[lb]{\smash{{\SetFigFont{9}{10.8}{\familydefault}{\color[rgb]{0,0,0}$X_3$}%
}}}}
\put(1646,-1193){\makebox(0,0)[lb]{\smash{{\SetFigFont{9}{10.8}{\familydefault}{\color[rgb]{0,0,0}$l_3$}%
}}}}
\put(2081,-706){\makebox(0,0)[lb]{\smash{{\SetFigFont{9}{10.8}{\familydefault}{\color[rgb]{0,0,0}$X_2$}%
}}}}
\put(1424,-981){\makebox(0,0)[lb]{\smash{{\SetFigFont{9}{10.8}{\familydefault}{\color[rgb]{0,0,0}$l_1$}%
}}}}
\put(1779,-820){\makebox(0,0)[lb]{\smash{{\SetFigFont{9}{10.8}{\familydefault}{\color[rgb]{0,0,0}$l_2$}%
}}}}
\end{picture}%
}
\end{center}
\caption{To definition of the bending energy for simplicial surfaces
\label{f.bending}}
\end{figure}
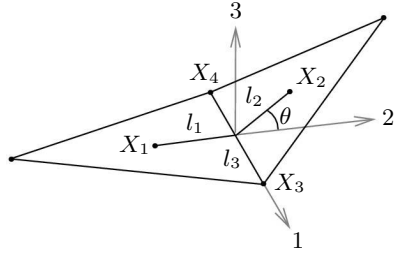

To derive the bending energy for simplicial surfaces let us
consider the limit of fine triangulations, where the angles
between the normals of neighboring triangles become small.
Consider an isometric deformation of two adjacent triangles. Let
$\theta$ be the external dihedral angle of the edge $e$, or,
equivalently, the angle between the normals of these triangles
(see Figure~\ref{f.bending}) and $\beta(\theta)$ the external
intersection angle between the circumcircles of the triangles (see
Figure~\ref{f.dEnergy}) as a function of $\theta$.

\begin{proposition}  \label{p.bending}
Assume that the circumcenters of two adjacent triangles do not
coincide. Then in the limit of small angles $\theta\to 0$ the
angle $\beta$ between the circles behaves as follows:
$$
\beta(\theta)=\beta(0)+\dfrac{l}{4L}\theta^2+o(\theta^3).
$$
Here $l$ is the length of the edge and $L\ne 0$ is the distance between the centers
of the circles.
\end{proposition}
\begin{proof}
Let us introduce the orthogonal coordinate system with the origin at the middle
point of the common edge $e$, the first basis vector directed along $e$, and the
third basis vector orthogonal to the left triangle. Denote by $X_1, X_2$ the
centers of the circumcircles of the triangles and by $X_3, X_4$ the end points of
the common edge (see Fig.~\ref{f.bending}). The coordinates of these points are
$X_1=(0,-l_1,0), X_2=(0,l_2\cos\theta,l_2\sin\theta), X_3=(l_3,0,0), ,
X_4=(-l_3,0,0)$. Here $2l_3$ is the length of the edge $e$, and $l_1$ and $l_2$ are
the distances from its middle point to the centers of the circumcirlces (for acute
triangles). The unit normals to the triangles are $N_1=(0,0,1)$ and
$N_2=(0,-\sin\theta,\cos\theta)$. The angle $\beta$ between the circumcircles
intersecting at the point $X_4$ is equal to the angle between the vectors
$A=N_1\times(X_4-X_1)$ and $B=N_2\times(X_4-X_2)$. The coordinates of these vectors
are $A=(-l_1,-l_3,0)$, $B=(l_2,-l_3\cos\theta,-l_3\sin\theta)$. This implies for
the angle
\begin{equation}            \label{beta-theta}
\cos\beta(\theta)=\dfrac{l_3^2\cos\theta-l_1l_2}{r_1r_2},
\end{equation}
where $r_i=\sqrt{l_i^2+l_3^2}, \ i=1,2$ are the radii of the
corresponding circumcircles. Thus $\beta(\theta)$ is an even
function, in particular
$\beta(\theta)=\beta(0)+B\theta^2+o(\theta^3)$. Differentiating
(\ref{beta-theta}) by $\theta^2$ we obtain
$$
B=\dfrac{l_3^2}{2r_1r_2\sin\beta(0)}.
$$
Also formula (\ref{beta-theta}) yields
$$
\sin\beta(0)=\dfrac{l_3L}{r_1r_2},
$$
where $L=| l_1+l_2 |$ is the distance between the centers of the circles. Finally
combining these formulas we obtain $B=l_3/(2L)$.
\end{proof}

This proposition motivates us to define the bending energy of
simplicial surfaces as
$$
E=\sum_{e\in E}\dfrac{l}{L}\theta^2.
$$
For discrete thin-shells this bending energy was suggested and analyzed by Grinspun
et al. \cite{GrinspunHDS,Grinspun}. The distance between the barycenters was used
for $L$ in the energy expression, and possible advantages in using circumcenters
were indicated. Numerical experiments demonstrate good qualitative simulation of
real processes.

 Further applications of the
discrete Willmore energy in particular for surface restoration, geometry denoising,
and smooth filling of a hole can be found in \cite{BobenkoSchroederWillmore}.

\section{Circular Nets as Discrete Curvature Lines}
\label{s.orthogonal}

Simplicial surfaces as studied in the previous section are too
unstructured for analytical investigation. An important tool in
the theory of smooth surfaces is the introduction of (special)
parametrizations of a surface. Natural analogues of parametrized
surfaces are quadrilateral surfaces, i.e. discrete surfaces made
from (not necessarily planar) quadrilaterals. The strips of
quadrilaterals obtained by gluing quadrilaterals along opposite
edges can be considered as coordinate lines on the quadrilateral
surface.

We start with a combinatorial description of the discrete surfaces under
consideration.

\begin{definition}
A cellular decomposition ${\mathcal D}$ of a two-dimensional
manifold (with boundary) is called a {\em quad-graph} if the cells
have four sides each.
\end{definition}
A quadrilateral surface is a mapping $f$ of a quad-graph to
${\mathbb R}^3$. The mapping $f$ is given just by the values at
the vertices of $\mathcal D$, and vertices, edges and faces of the
quad-graph and of the quadrilateral surface correspond.
Quadrilateral surfaces with planar faces were suggested by Sauer
\cite{Sauer} as discrete analogs of conjugate nets on smooth
surfaces. The latter are the mappings $(x,y)\mapsto
f(x,y)\in{\mathbb R}^3$ such that the mixed derivative $f_{xy}$ is
tangent to the surface.

\begin{definition}\label{d.circular}
A quadrilateral surface $f:{\mathcal D}\to {\mathbb R}^3$ all
faces of which are circular (i.e. the four vertices of each face
lie on a common circle) is called a {\em circular net} (or
discrete orthogonal net).
\end{definition}

Circular nets as discrete analogues of curvature line parametirzed
surfaces were mentioned by Martin, de Pont, Sharrock and Nutbourne
\cite{MartinPS, Nutbourne} . The curvature lines on smooth
surfaces continue through any point. Keeping in mind the analogy
to the curvature line parametrized surfaces one may in addition
require that all vertices of a circular net are of even degree.

\begin{figure}[h]
\begin{center}
\parbox[c]{0.45\textwidth}{{\includegraphics[width=0.45\textwidth]{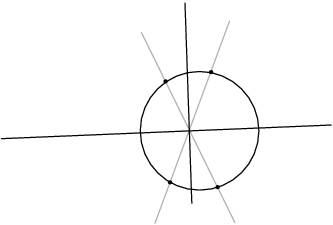}}}
\end{center}
\caption{Principal curvature directions of a circular quadrilateral
\label{f.orthogonal}}
\end{figure}

A smooth conjugate net $f:D\to {\mathbb R}^3$ is a curvature line
parametrization if and only if it is orthogonal. The bisectors of
the diagonals of a circular quadrilateral intersect orthogonally
(see Figure~\ref{f.orthogonal}) and can be interpreted
\cite{BobenkoTsarev} as discrete principal curvature directions.

There are deep reasons to treat circular nets as a discrete
curvature line parametrization.
\begin{itemize}
    \item The class of circular nets as well as the class of curvature line
    parametrized surfaces is invariant under M\"obius transformations.
    \item Take an infinitesimal quadrilateral $(f(x,y), f(x+\epsilon), y),
    f(x+\epsilon),y+\epsilon), f(x, y+\epsilon))$ of a curvature line parametrized
    surface. A direct computation (see \cite{BobenkoTsarev}) shows that in the limit $\epsilon\to 0$
    the imaginary part of its cross-ratio is of order $\epsilon^3$. Note that circular
    quadrilaterals are characterized by having real cross-ratios.
    \item For any smooth curvature line parametrized surface
    $f:D\to {\mathbb R}^3$ there exists a family of discrete circular nets
    converging to $f$. Moreover the convergence is $C^\infty$, i.e. with all
    derivatives. The details can be found in \cite{BobenkoMatthesSuris}.
\end{itemize}

One more argument in favor of Definition \ref{d.circular} is that circular nets
satisfy the {\em consistency principle}, which has proven to be one of the
organizing principles in discrete differential geometry \cite{BobenkoSurisDDG}. The
consistency principle singles out fundamental geometries by the requirement that
the geometry can be consistently extended to a combinatorial grid one dimension
higher. The consistency of circular nets was shown by Doliwa and Santini
\cite{DoliwaSantini} based on the following classical theorem.

\begin{figure}[h]
\begin{center}
\parbox[c]{0.35\textwidth}{{\includegraphics[width=0.35\textwidth]{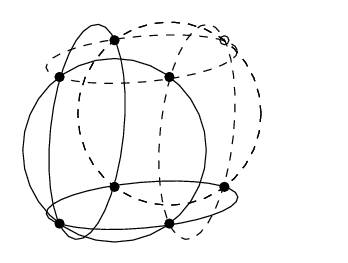}}}
\hspace{1.5cm}
\parbox[c]{0.25\textwidth}{{\includegraphics[width=0.25\textwidth]{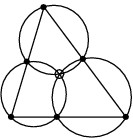}}}
\end{center}
\caption{Miquel's theorem: spherical and planar versions\label{f.Miquel}}
\end{figure}

\begin{theorem}(Miquel)
Consider a combinatorial cube in ${\mathbb R}^3$ with planar
faces. Assume that three neighboring faces of the cube are
circular. Then the remaining three faces are also circular.
\end{theorem}
Equivalently, provided the four-tuples of black vertices
corresponding to three neighboring faces of the cube lie on
circles, the three circles determined by the triples of points
corresponding to three remaining faces of the cube all intersect
(at the white vertex in Figure~\ref{f.Miquel}). It is easy to see
that all vertices of Miquel's cube lie on a sphere. Mapping the
vertex shared by the three original faces to infinity by a
M\"obius transformation, we obtain an equivalent planar version of
Miquel's theorem. This version, also shown in
Figure~\ref{f.Miquel}, can be proven by means of elementary
geometry.

Finally note that circular nets are also treated as a
discretization of triply orthogonal coordinate systems. Triply
orthogonal coordinate systems in ${\mathbb R}^3$ are maps
$(x,y,z)\mapsto f(x,y,z)\in{\mathbb R}^3$ from a subset of
${\mathbb R}^3$ with mutually orthogonal $f_x, f_y, f_z$. Due to
the classical Dupin theorem, the level surfaces of a triply
orthogonal coordinate system intersect along their common
curvature lines. Accordingly, discrete triply orthogonal systems
are defined as maps from ${\mathbb Z}^3$ (or a subset thereof) to
${\mathbb R}^3$ with all elementary hexahedra lying on spheres
\cite{Bobenko99}. Due to Miquel's theorem a discrete orthogonal
system is uniquely determined by the circular nets corresponding
to its coordinate two-planes (see \cite{DoliwaSantini} and
\cite{BobenkoSurisDDG}).

\section{Discrete Isothermic Surfaces}
\label{s.isothermic}

In this section and in the following one, we investigate discrete
analogs of special classes of surfaces obtained by imposing
additional conditions in terms of circles and spheres.

\begin{figure}[ht]
\begin{center}
\parbox[c]{0.2\textwidth}{{\includegraphics[width=0.2\textwidth]{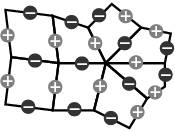}}}
\end{center}
\caption{Labelling the edges of a discrete isothermic surface\label{f.labelling}}
\end{figure}

\begin{figure}[ht]
\begin{center}
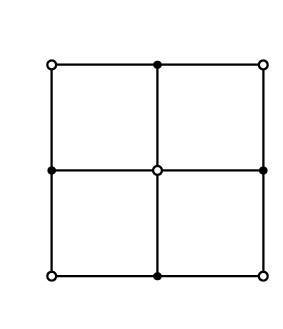\hspace{1cm}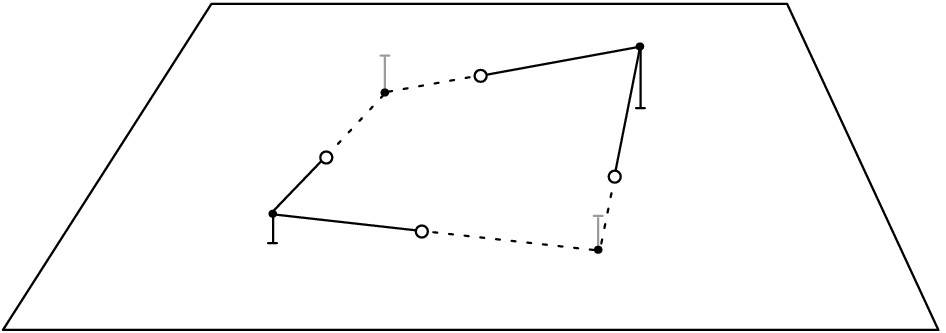
\end{center}
\caption{Central spheres of a discrete isothermic surface: combinatorics
{\em(left)}, and the M\"obius normalized picture  for $N=2$ {\em(right)}.
\label{f.central-sphere}}
\end{figure}

We start with minor combinatorial restrictions. Suppose that the
vertices of $\mathcal D$ are colored black or white so that the
two ends of each edge have different colors. Such a coloring is
always possible for topological discs. To model the curvature
lines, suppose also that the edges of a quad-graph $\mathcal D$
may consistently be labelled `$+$' and `$-$', as in
Figure~\ref{f.labelling} (For this it is necessary that each
vertex has an even number of edges).
Let $f_0$ be a vertex of a circular net, $f_1,f_3,\ldots,f_{4N-1}$
be its neighbors, and $f_2,f_4,\ldots,f_{4N}$ its next-neighbors
(see Figure~\ref{f.central-sphere}~(left)). We call the vertex
$f_0$ generic if it is not co-spherical with all its neighbors and
a circular net $f:{\mathcal D}\to{\mathbb R}^3$ generic if all its
vertices are generic.

Let  $f:{\mathcal D}\to{\mathbb R}^3$ be a generic circular net
such that every its vertex is co-spherical with all its
next-neighbors. We will call the corresponding sphere {\em
central}. For analytical description of this geometry let us map
the vertex $f_0$ to infinity by a M\"obius transformation
${\mathcal M}(f_0)=\infty$, and denote by $F_i={\mathcal M}(f_i),
i=1,\ldots,4N$ the M\"obius images of $f_i$. The points
$F_2,F_4,\ldots,F_{4N}$ are obviously coplanar. The circles of the
faces are mapped to straight lines. For the cross-ratios we get
$$
q(f_0,f_{2k-1},f_{2k},f_{2k+1})=\dfrac{F_{2k}-F_{2k+1}}{F_{2k}-F_{2k-1}}=
\dfrac{z_{2k+1}}{z_{2k-1}},
$$
where $z_{2k+1}$ is the coordinate of $F_{2k+1}$ orthogonal to the
plane $\mathcal P$ of $F_2,F_4,\ldots,F_{4N}$.  (Note that since
$f_0$ is generic none of the $z_i$ vanishes). As a corollary we
get for the product of all cross-ratios:
\begin{equation}\label{product_isothermic}
\prod_{k=1}^n q(f_0,f_{2k-1},f_{2k},f_{2k+1})=1.
\end{equation}
\begin{definition}
A circular net $f:{\mathcal D}\to{\mathbb R}^3$ satisfying condition
(\ref{product_isothermic}) at each vertex is called a {\em discrete isothermic
surface}.
\end{definition}
This definition was first suggested in \cite{BobenkoPinkall-Isoth}
for the case of the combinatorial square grid ${\mathcal
D}={\mathbb Z}^2$. In this case if the vertices are labelled by
$f_{m,n}$ and the corresponding cross-ratios by $q_{m,n}:=
q(f_{m,n},f_{m+1,n},f_{m+1,n+1},f_{m,n+1})$, the condition
(\ref{product_isothermic}) reads
$$
q_{m,n}q_{m+1,n+1}=q_{m+1,n}q_{m,n+1}.
$$

\begin{proposition}
Every vertex $f_{m,n}$ of a discrete isothermic surface
$f:{\mathbb Z}^2\to{\mathbb R}^3$ possesses a central sphere, i.e.
the points $f_{m,n},f_{m-1,n-1},f_{m+1,n-1},f_{m+1,n+1}$ and
$f_{m-1,n+1}$ are co-spherical. Moreover for generic circular maps
$f:{\mathbb Z}^2\to{\mathbb R}^3$ this property characterizes
discrete isothermic surfaces.
\end{proposition}

\begin{proof} Use the notation of Figure~\ref{f.central-sphere}, with $f_0\equiv f_{m,n}$,
and the same argument with the M\"obius transformation $\mathcal
M$ which maps $f_0$ to $\infty$. Consider the plane $\mathcal P$
determined by the points $F_2, F_4$ and $F_6$. Let as above $z_k$
be the coordinates of $F_k$ orthogonal to the plane $\mathcal P$.
Condition (\ref{product_isothermic}) yields
$$
\dfrac{F_8-F_1}{F_8-F_7}= \dfrac{z_1}{z_7}.
$$
This implies that the $z$-coordinate of the point $F_8$ vanishes, thus
$F_8\in{\mathcal P}$.
\end{proof}

 The property to be discrete isothermic is also 3D-consistent, i.e. can be
consistently imposed on all faces of a cube. This was proven first by
Hertrich-Jeromin, Hoffmann and Pinkall \cite{JerominHoffmannPinkall} (see also
\cite{BobenkoSurisDDG} for generalizations and modern treatment).

An important subclass of discrete isothermic surfaces is given by
the condition that all the faces are conformal squares, i.e. their
cross ratio equals $-1$. All conformal squares are M\"obius
equivalent, in particular equivalent to the standard square. This
is a direct discretization of the definition of smooth isothermic
surfaces. The latter are immersions $(x,y)\mapsto
f(x,y)\in{\mathbb R}^3$ satisfying
\begin{equation}            \label{isoth-smooth}
\| f_x\|=\| f_y\|, f_x\perp f_y, f_{xy}\in {\rm span}\{f_x,f_x\},
\end{equation}
i.e. conformal curvature line parametrizations. Geometrically this definition means
that the curvature lines divide the surface into infinitesimal squares.
\begin{figure}[ht]
\begin{center}
\parbox[c]{0.2\textwidth}{{\includegraphics[width=0.2\textwidth]{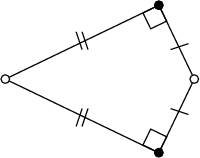}}}
\end{center}
\caption{Right-angled kites are conformal squares \label{f.kite}}
\end{figure}

\begin{figure}
 \begin{center}
\parbox[c]{0.25\textwidth}{\includegraphics[width=0.25\textwidth]{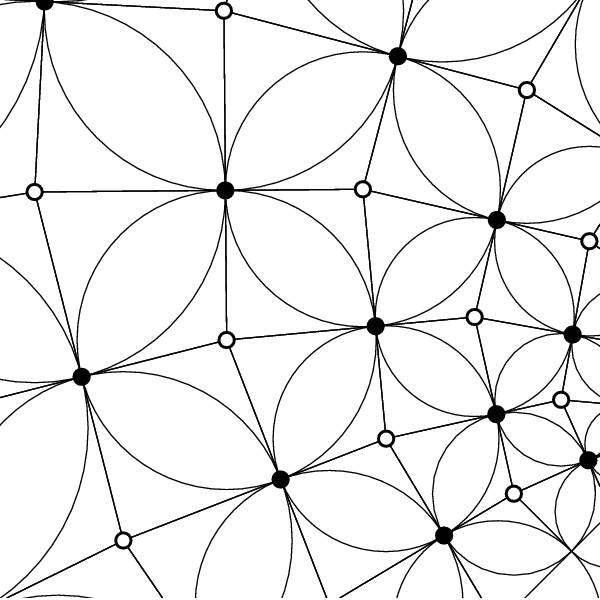}}
 \hspace{2cm}
\parbox[c]{0.35\textwidth}{\includegraphics[width=0.3\textwidth]{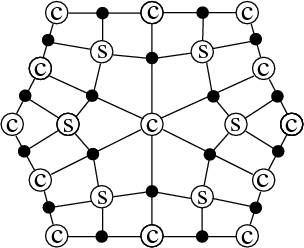}}
   \end{center}
 \caption{To the definition of discrete S-isothermic surfaces:
 orthogonal circle patterns as discrete conformal maps {\em(left)} and
 combinatorics of S-quad-graphs {\em(right)}.}
 \label{f.S-quad-graph}
\end{figure}

The class of discrete isothermic surfaces is too general and the
surfaces are not rigid enough. In particular one can show that the
surface can vary preserving all its black vertices. In this case,
one white vertex can be chosen arbitrarily
\cite{BobenkoPinkall-DIGP}. Thus, we introduce a more rigid
subclass. To motivate its definition, let us look at the problem
of discretizing the class of conformal maps $f:D\to{\mathbb C}$
for $D\subset{\mathbb C}={\mathbb R}^2$. Conformal maps are
characterized by the conditions
\begin{equation}
\label{eq:conformality conditions}
 |f_x|=|f_y|,\qquad f_x\perp f_y.
\end{equation}
To define discrete conformal maps $f:\mathbb{Z}^2\supset D\to \mathbb{C}$, it is
natural to impose these two conditions on two different sub-lattices (white and
black) of $\mathbb{Z}^2$, i.e.\ to require that the edges meeting at a white vertex
have equal length and the edges at a black vertex meet orthogonally. This
discretization leads to the circle patterns with the combinatorics of the square
grid introduced by Schramm \cite{Schramm-SG}. Each circle intersects four
neighboring circles orthogonally and the neighboring circles touch cyclically
(Fig.\ref{f.S-quad-graph} {\em(left)}).

The same properties imposed for quadrilateral surfaces with the
combinatorics of the square grid $f:\mathbb{Z}^2\supset{\mathcal
D}\to {\mathbb R}^3$ lead to an important subclass of discrete
isothermic surfaces. Let us require for a discrete quadrilateral
surface that:
\begin{itemize}
    \item the faces are orthogonal kites,
    \item the edges meet at black vertices orthogonally
    (black vertices are at orthogonal corners of the kites),
    \item the kites which do not share a common vertex are not coplanar (locality condition).
\end{itemize}

Observe that the orthogonality condition (at black vertices)
implies that one pair of opposite edges meeting at a black vertex
lies on a straight line. Together with the locality condition this
implies that there are two kinds of white vertices, which we
denote by $\cir$ and $\sph$. Each kite has white vertices of both
types and the kites sharing a white vertex of the first kind
$\cir$ are coplanar.

These conditions imposed on the quad-graphs lead to S-quad-graphs
and S-isothermic surfaces (the latter were introduced in
\cite{BobenkoPinkall-DIGP} for the combinatorics of the square
grid).

\begin{definition}
\label{d.S-isothermic} An {\em S-quad-graph} $\mathcal D$ is a
quad-graph with black and two kinds of white vertices such that
the two ends of each edge have different colors and each
quadrilateral has vertices of all kinds (see
Fig.~\ref{f.S-quad-graph} {\em(right)}). Let $V_{b}({\mathcal D})$
be the set of black vertices. A {\em discrete S-isothermic
surface} is a map
\begin{equation*}
  f_{b}:V_{b}({\mathcal D})\rightarrow\mathbb{R}^3,
\end{equation*}
with the following properties:
\begin{enumerate}
\renewcommand{\theenumi}{{\em(\roman{enumi})}}
\renewcommand{\labelenumi}{\theenumi}
\item If $v_1, \ldots, v_{2n}\in V_{b}({\mathcal D})$ are the neighbors of a
  $\cir$-labeled vertex in cyclic order, then $f_{b}(v_1), \ldots,
  f_{b}(v_{2n})$ lie on a circle in $\mathbb{R}^3$ in the same cyclic order.
  This defines a map from the $\cir$-labeled vertices to the set of circles
  in $\mathbb{R}^3$.
\item If $v_1, \ldots, v_{2n}\in V_{b}({\mathcal D})$ are the neighbors of an
  $\sph$-labeled vertex, then $f_{b}(v_1),\ldots,$ $f_{b}(v_{2n})$ lie on a
  sphere in $\mathbb{R}^3$. This defines a map from the $\sph$-labeled
  vertices to the set of spheres in $\mathbb{R}^3$.
\item If $v_c$ and $v_s$ are the $\cir$-labeled and the $\sph$-labeled vertex
  of a quadrilateral of $\mathcal D$, then the circle corresponding to $v_c$
  intersects the sphere corresponding to $v_s$ orthogonally.
\end{enumerate}
\end{definition}
Discrete S-isothermic surfaces are therefore composed of tangent
spheres and tangent circles, with the spheres and circles
intersecting orthogonally. The class of discrete S-isothermic
surfaces is obviously invariant under M\"obius transformations.

Given a discrete S-isothermic surface, one can add the centers of
the spheres and circles to it giving a map $V({\mathcal
D})\rightarrow\mathbb{R}^3$. The discrete isothermic surface
obtained is called the {\em central extension} of the discrete
S-isothermic surface. All its faces are orthogonal kites.

An important fact of the theory of isothermic surfaces (smooth and discrete) is the
existence of a dual isothermic surface \cite{BobenkoPinkall-Isoth}. Let $f:{\mathbb
R}^2\supset D\to {\mathbb R}^3$ be an isothermic immersion. Then the formulas
$$
f^*_x=\dfrac{f_x}{\|f_x\|^2},\quad f^*_y=-\dfrac{f_y}{\|f_y\|^2}.
$$
define an isothermic immersion $f^*:{\mathbb R}^2\supset D\to {\mathbb R}^3$ which
is called the {\em dual isothermic surface}. Indeed, one can easily check that the
form $df^*$ is closed and $f^*$ satisfies (\ref{isoth-smooth}). Exactly the same
formulas can be applied in the discrete case.

\begin{proposition}
Suppose $f:{\mathcal D}\to {\mathbb R}^3$ is a discrete isothermic
surface and suppose the edges have been consistently labeled `$+$'
and `$-$', as in Fig.~\ref{f.labelling}. Then the {\em dual
discrete isothermic surface} $f^*:{\mathcal D}\to {\mathbb R}^3$
is defined by the formula
$$
\Delta f^*=\pm \dfrac{\Delta f}{\|\Delta f\|^2},
$$
where $\Delta f$ denotes the difference of values at neighboring
vertices and the sign is chosen according to the edge label.
\end{proposition}
The closedness of the corresponding discrete form is elementary
checked for one kite.

\begin{proposition}
The dual of the central extension of a discrete S-isothermic
surface is the central extension of another discrete S-isothermic
surface.
\end{proposition}
If we disregard the centers, we obtain the definition of the {\em
dual discrete S-isothermic surface}. The dual discrete
S-isothermic surface can be defined also without referring to the
central extension \cite{BobenkoHoffmannSpringborn}.

\section{Discrete Minimal Surfaces and Circle Patterns: Geometry from Combinatorics}
\label{s.minimal}

In this section (following \cite{BobenkoHoffmannSpringborn}) we
define {\em discrete minimal S-isothermic surfaces} (or {\em
discrete minimal surfaces} for short) and present the most
important facts from their theory. The main idea of the approach
of \cite{BobenkoHoffmannSpringborn} is the following. Minimal
surfaces are characterized among (smooth) isothermic surfaces by
the property that they are dual to their Gauss map. The duality
transformation and this characterization of minimal surfaces
carries over to the discrete domain. The role of the Gauss map is
played by discrete S-isothermic surfaces whose spheres all
intersect one fixed sphere orthogonally.

\subsection{Koebe polyhedra}

A {\em circle packing} in $S^2$ is a configuration of disjoint discs which may
touch but not intersect. The construction of discrete S-isothermic ``round
spheres'' is based on their relation to circle packings in $S^2$.  The following
theorem is central in this theory.

\begin{theorem}
\label{t.orthoKoebe} For every polytopal\footnote{We call a cellular
  decomposition of a surface {\em polytopal}, if the closed cells are closed
  discs, and two closed cells intersect in one closed cell if at all.}
cellular decomposition of the sphere, there exists a pattern of circles in the
sphere with the following properties.  There is a circle corresponding to each face
and to each vertex.  The vertex circles form a packing with two circles touching if
and only if the corresponding vertices are adjacent. Likewise, the face circles
form a packing with circles touching if and only if the corresponding faces are
adjacent. For each edge, there is a pair of touching vertex circles and a pair of
touching face circles. These pairs touch in the same point, intersecting each other
orthogonally.

This circle pattern is unique up to M{\"o}bius transformations.
\end{theorem}

The first published statement and proof of this theorem is contained in
\cite{BrightwellScheinerman}. For generalizations, see \cite{Schramm-Egg},
\cite{Rivin}, and \cite{BobenkoSpringborn-Variational}, the latter also for a
variational proof.

Theorem~\ref{t.orthoKoebe} is a generalization of the following remarkable
statement about circle packings due to Koebe \cite{Koebe}.
\begin{theorem}{\bf (Koebe)}
\label{t.koebe}
 For every triangulation of the sphere there is a
packing of circles in the sphere such that circles correspond to vertices, and two
circles touch if and only if the corresponding vertices are adjacent. This circle
packing is unique up to M{\"o}bius transformations of the sphere.
\end{theorem}

\begin{figure}[tbp]
\hfill%
\includegraphics[width=0.3\textwidth]{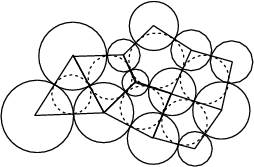}%
\hfill%
\includegraphics[width=0.3\textwidth]{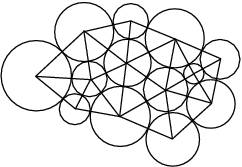}%
\hfill%
\includegraphics[width=0.3\textwidth]{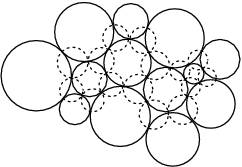}%
\hspace*{\fill}%
\caption{{\em Left:} An orthogonal circle pattern
    corresponding to a cellular decomposition.
    {\em Middle:} A circle packing corresponding to a triangulation.
    {\em Right:}  The orthogonal circles. }
 \label{f.Koebe-combinatorics}
\end{figure}

Observe that, for a triangulation, one automatically obtains not one but two
orthogonally intersecting circle packings, as shown in
Fig.~\ref{f.Koebe-combinatorics} {\em(right)}. Indeed, the circles passing through
the points of contact of three mutually touching circles intersect these
orthogonally, and thus Theorem~\ref{t.koebe} is a special case of
Theorem~\ref{t.orthoKoebe}.

Consider a circle pattern of Theorem~\ref{t.orthoKoebe}.
Associating white vertices to circles and black vertices to their
intersection points, one obtains a quad-graph. Actually we have an
S-quad-graph: Since the circle pattern is comprised by two circle
packings intersecting orthogonally we have two kinds of white
vertices, $\sph$ and $\cir$, corresponding to the circles of the
two packings.

Now let us construct the spheres intersecting $S^2$ orthogonally
along the circles marked by $\sph$. If we then connect the centers
of touching spheres, we obtain a convex polyhedron, all of whose
edges are tangent to the sphere $S^2$. Moreover, the circles
marked with $\cir$ are inscribed in the faces of the polyhedron.
Thus we have a discrete S-isothermic surface.

We arrive at the following theorem, which is equivalent to
Theorem~\ref{t.orthoKoebe} (see also \cite{Ziegler}).
\begin{theorem}
\label{t.polyKoebe}
Every polytopal cell decomposition of the sphere can be
realized by a polyhedron with edges tangent to the sphere. This realization is
unique up to projective transformations which fix the sphere.
\end{theorem}
These polyhedra are called the {\em Koebe polyhedra}. We interpret the
corresponding discrete S-isothermic surfaces as conformal discretizations of the
``round sphere''.

\subsection{Definition of discrete minimal surfaces}

\begin{figure}[htb]
  \centering
  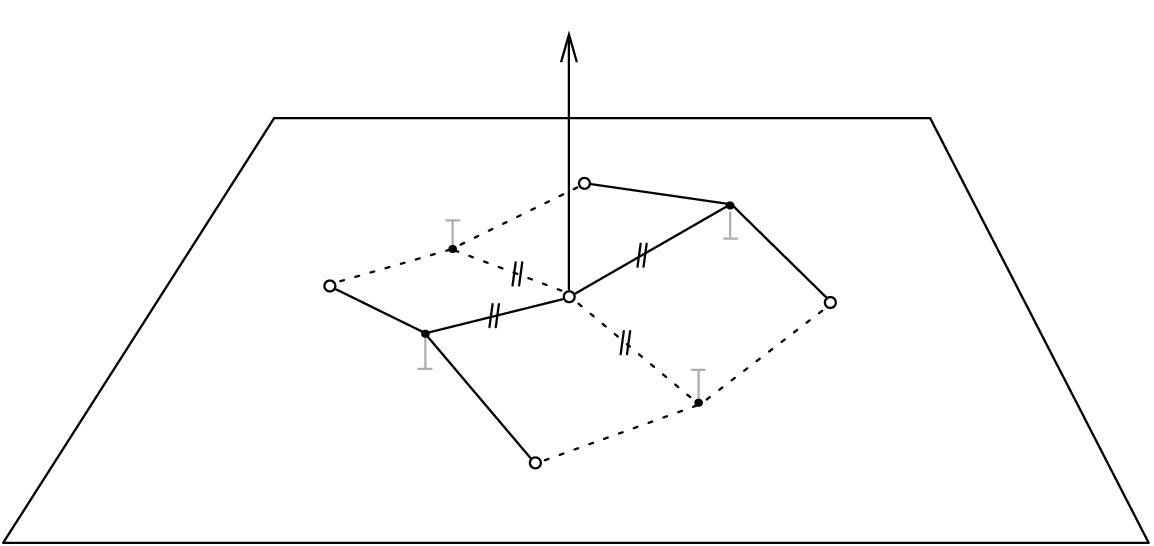
  \caption{To the definition of discrete minimal surfaces: the
  tangent plane through the center $f(x)$ of a sphere and the centers
  $f(x_j)$ of the neighboring circles. The circles and the sphere
  intersect orthogonally at the black points $f(y_j)$.}
  \label{f.discrete_minimal_definition}
\end{figure}

Let $f:{\mathcal D}\to {\mathbb R}^3$ be a discrete S-isothermic
surface. Suppose $x\in {\mathcal D}$ is a white $\sph$ vertex of
the quad-graph $\mathcal D$, i.e. $f(x)$ is the center of a
sphere. Consider all quadrilaterals of $\mathcal D$ incident to
$x$ and denote by $y_1\ldots y_{2n}$ their black vertices and by
$x_1\ldots x_{2n}$ their white $\cir$ vertices. We will call the
vertices $x_1\ldots x_{2n}$ neghboring $x$ in $\mathcal D$.
(Generically, $n=2$.) Then $f(y_j)$ are the points of contact with
the neighboring spheres and simultaneously points of intersection
with the orthogonal circles centered at $f(x_j)$ (see
Fig.~\ref{f.discrete_minimal_definition}).

Consider only the white vertices of the quad-graph. Observe that
each $\cir$ vertex of a discrete S-isothermic surface and all its
$\sph$ neighbors are coplanar. Indeed the plane of the circle
centered at $f(\cir)$ contains all its $f(\sph)$ neighbors. The
same condition imposed at the $\sph$ vertices leads to a special
class of surfaces.

\begin{definition}
  \label{d.discrete_minimal}
 A discrete S-isothermic surface $f:{\mathcal D}\to {\mathbb R}^3$ is called
 {\em discrete minimal} if for each sphere of
 $f({\mathcal D})$ there exists a plane which contains the center
  $f(\sph)$ of the sphere as well as the centers $f(\cir)$ of all neighboring
 circles orthogonal to the sphere, i.e. if each white vertex of $f({\mathcal D})$
 is coplanar to all its white neighbors. These planes should be
 considered as tangent planes to the discrete minimal surface.
\end{definition}

\begin{theorem}               \label{t.minimal_dual_koebe}
  An S-isothermic discrete surface $f$ is a discrete minimal surface, if and only
  if the dual S-isothermic surface $f^*$ corresponds to a Koebe polyhedron. In
  that case the dual surface $N:=f^*:V_w({\mathcal D})\to {\mathbb R}^3$ at white
  vertices $V_w({\mathcal D})$ can be treated as the Gauss map $N$ of
  the discrete minimal surface: At \sph vertices $N$ is orthogonal
  to the tangent planes, and at \cir vertices $N$ is orthogonal to the planes of
  the circles centered at $f(\cir)$.
\end{theorem}

\begin{proof}
That the S-isothermic dual of a Koebe polyhedron is a discrete minimal
  surface is fairly obvious. On the other hand, let
  $f:{\mathcal D}\rightarrow {\mathbb R}^3$ be
  a discrete minimal surface and let $x\in {\mathcal D}$ and
  $y_1\ldots y_{2n}\in {\mathcal D}$ be
  as above.  Let
  $f^*:{\mathcal D}\rightarrow {\mathbb R}^3$ be the dual
  S-isothermic surface. We need
  to show that all circles of $f^*$ lie in one and the same sphere
  $S$ and that all the spheres of $f^*$ intersect $S$ orthogonally.
  Since the quadrilaterals of a discrete S-isothermic surface are kites the
minimality condition of Definition~\ref{d.discrete_minimal} can be reformulated as
follows: There is  a plane through $f(x)$ such that the points $\{f(y_j)\;|\;j
\text{ even}\}$ and the points $\{f(y_j)\;|\;j \text{ odd}\}$ lie in planes which
are parallel to it at the same distance on opposite sides (see
Fig.~\ref{f.discrete_minimal_definition}). It follows immediately that the
  points $f^*(y_1)\ldots f^*(y_{2n})$ lie on a circle
  $c_x$ in a sphere $S_x$ around $f^*(x)$. The plane of $c_x$ is orthogonal
  to the normal $N$ to the tangent plane at $f(x)$.  Let $S$ be the sphere
  which intersects $S_x$ orthogonally in $c_x$. The orthogonal circles
  through $f^*(y_1)\ldots f^*(y_{2n})$ also lie in $S$.
  Hence, all spheres of $f^*$ intersect $S$ orthogonally and all
  circles of $f^*$ lie in $S$.
\end{proof}

Theorem~\ref{t.minimal_dual_koebe} is a complete analogue of Christoffel's
characterization \cite{Christoffel} of continuous minimal surfaces.

\begin{theorem}{\bf (Christoffel)}
  Minimal surfaces are isothermic. An isothermic immersion is a minimal
  surface, if and and only if the dual immersion is contained in a sphere. In
  that case the dual immersion is in fact the Gauss map of the minimal
  surface.
\end{theorem}

Thus a discrete minimal surface is a discrete S-isothermic surface
which is dual to a Koebe polyhedron; the latter is its Gauss map
and is a discrete analogue of a conformal parametrization of the
sphere.

The simplest infinite orthogonal circle pattern in the plane
consists of circles with equal radius $r$ and centers on a square
grid with spacing $\frac{1}{2}\sqrt{2}\,r$. One can project it
stereographically to the sphere, construct orthogonal spheres
through half of the circles and dualize to obtain a discrete
version of Enneper's surface. See
Fig.~\ref{f.dSurfaces}{\em(right)}. Only the circles are shown.

\subsection{Construction of discrete minimal surfaces}

\begin{figure}[htb]
  \centering
  \parbox[b]{0.3\textwidth}{\includegraphics[width=0.3\textwidth]{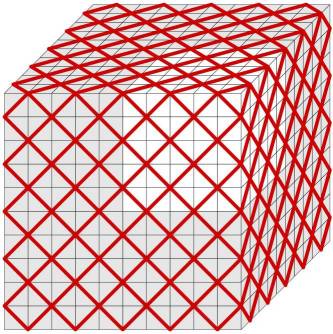}\\
    \centerline{Gauss image of the curvature lines}}
   \raisebox{0.2\textwidth}{\huge\ $\rightarrow{}$ }
   \parbox[b]{0.3\textwidth}{\includegraphics[width=0.3\textwidth]{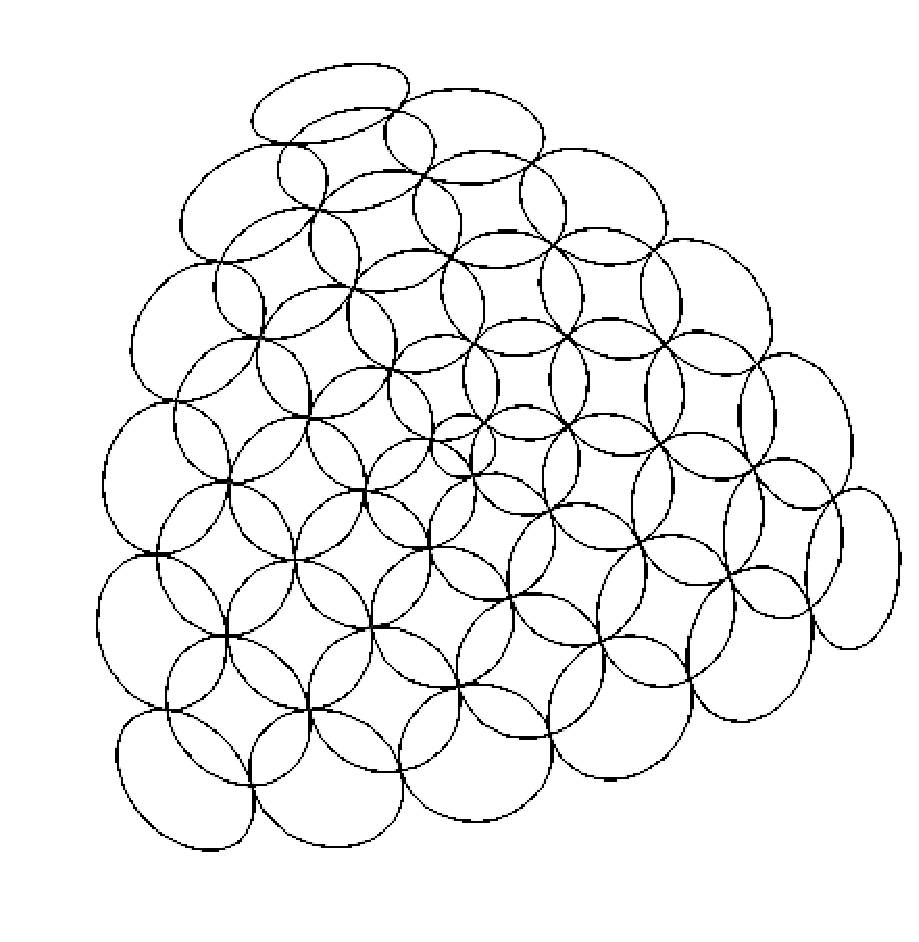}\\
     \centerline{circle pattern}}  \\
   \vspace{\baselineskip}
   \raisebox{0.15\textwidth}{\huge\ $\rightarrow{}$ }
   \parbox[b]{0.25\textwidth}{\includegraphics[width=0.25\textwidth]{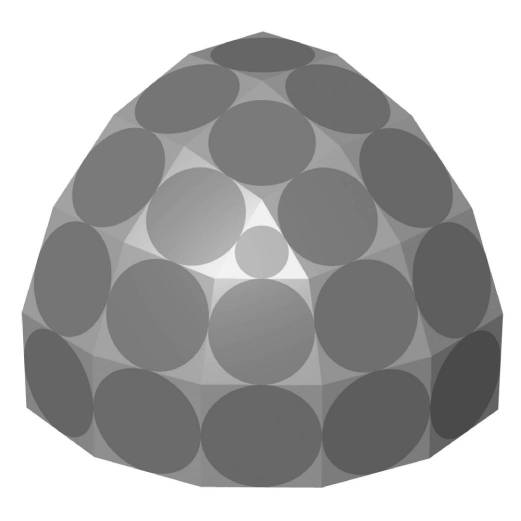}\\
     \centerline{Koebe polyhedron}}
   \raisebox{0.15\textwidth}{\huge\ $\rightarrow{}$ }
   \parbox[b]{0.35\textwidth}{\raisebox{0.03\textwidth}{\includegraphics[width=0.35\textwidth]{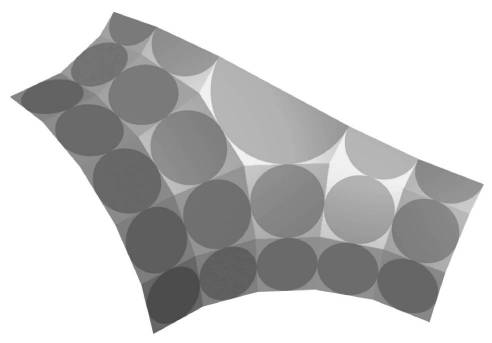}}\\
   \centerline{discrete minimal surface}}
  \caption{Construction of the discrete Schwarz P-surface.}
  \label{f.recipe}
\end{figure}

\begin{figure}[htb]
\centering%
\includegraphics[width=0.57\textwidth,clip=true]{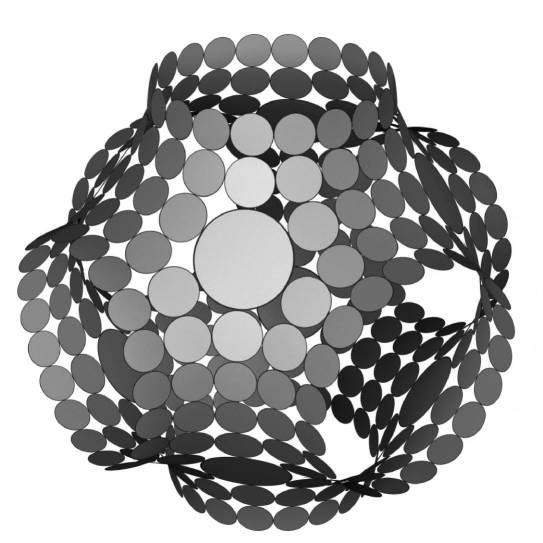}%
\caption{A discrete minimal Schwarz-P surface.} \label{f.schwarz}
\end{figure}

A general method to construct discrete minimal surfaces is
schematically shown in the following diagram.
\begin{center}
  continuous minimal surface \\
  {\large$\Downarrow$} \\
  image of curvature lines under Gauss-map\\
  {\large$\Downarrow$} \\
  cell decomposition of (a branched cover of) the sphere\\
  {\large$\Downarrow$} \\
  orthogonal circle pattern \\
  {\large$\Downarrow$} \\
  Koebe polyhedron\\
  {\large$\Downarrow$} \\
  discrete minimal surface
\end{center}
As usual in the theory on minimal surfaces \cite{HoffmannKarcher},
one starts constructing such a surface with a rough idea of how it
should look. To use our method, one should understand its Gauss
map and the {\em combinatorics}\/ of the curvature-line pattern.
The image of the curvature-line pattern under the Gauss map
provides us with a cell decomposition of (a part of) $S^2$ or a
covering.  From these data, applying Theorem~\ref{t.orthoKoebe},
we obtain a Koebe polyhedron with the prescribed combinatorics.
Finally, the dualization step yields the desired discrete minimal
surface. For the discrete Schwarz P-surface the construction
method is demonstrated in Fig.~\ref{f.recipe} and
Fig.~\ref{f.schwarz}.

Let us emphasize that our data, aside from possible boundary
conditions, are purely {\em combinatorial}---the combinatorics of
the curvature line pattern. All faces are quadrilaterals and
typical vertices have four edges. There may exist distinguished
vertices (corresponding to the ends or umbilic points of a minimal
surface) with a different number of edges.

The most nontrivial step in the above construction is the third
one listed in the diagram. It is based on the generalized Koebe
Theorem~\ref{t.orthoKoebe}. It implies the existence and
uniqueness for the discrete minimal S-isothermic surface under
consideration, but not only this. A constructive proof of the
generalized Koebe theorem suggested in
\cite{BobenkoSpringborn-Variational} is based on a variational
principle and also provides a method for the numerical
construction of circle patterns. An alternative algorithm by
Thurston was implemented in Stephenson's program {\tt{}circlepack}
(see \cite{Stephenson} for an exhaustive presentation of the
theory of circle packings and their numerics). The first step is
to transfer the problem from the sphere to the plane by a
stereographic projection. Then the radii of the circles are
calculated. If the radii are known, it is easy to reconstruct the
circle pattern. The above-mentioned variational method is based on
the observation that the equations for the radii are the equations
for a critical point of a convex function of the radii. The
variational method involves minimizing this function to solve the
equations.

Let us describe the variational method of \cite{BobenkoSpringborn-Variational} for
construction of (orthogonal) circle patterns in the plane and demonstrate how it
can be applied to construct the discrete Schwarz-P surface.

Instead of the radii $r$ of the circles, we use the logarithmic radii
\[ \rho=\log r. \]
For each circle $j$, we need to find a $\rho_j$ such that the corresponding radii
solve the circle pattern problem. This leads to the following equations, one for
each circle. The equation for circle $j$ is
\begin{equation}
  \label{closure}
  2\sum_{\makebox[0pt]{\scriptsize$\text{neighbors } k$}}
  \arctan e^{\rho_k-\rho_j} = \Phi_j,
\end{equation}
where the sum is taken over all neighboring circles $k$. For each circle $j$,
$\Phi_j$ is the nominal angle covered by the neighboring circles. It is normally
$2\pi$ for interior circles, but it differs for circles on the boundary or for
circles where the pattern branches.

\begin{theorem}
The critical points of the functional
  \begin{eqnarray*}
    S(\rho) = \sum_{(j,k)}\big(
    \im\Li_2(ie^{\rho_k-\rho_j}) +
    \im\Li_2(ie^{\rho_j-\rho_k}) - {\pi\over 2}(\rho_j+\rho_k)
    \big)
    + \sum_j \Phi_j\rho_j.
  \end{eqnarray*}
  correspond to orthogonal circle patterns in the plane with cone angles $\Phi_j$ at
  the centers of circles ($\Phi_j=2\pi$ for internal circles).
  Here, the first sum is taken over all pairs $(j,k)$ of neighboring circles,
  the second sum is taken over all circles $j$, and the dilogarithm function
  $\Li_2(z)$ is defined by $\Li_2(z)=-\int_0^z \log(1-\zeta)\,d\zeta/\zeta$.
  The functional is scale-invariant and restricted to the subspace $\sum_j \rho_j=0$
it is convex.
\end{theorem}
\begin{proof} The formula for the functional follows from (\ref{closure}) and
$$
\dfrac{d}{dx}\im\Li_2(ie^x)=\dfrac{1}{2i}\log\dfrac{1+ie^x}{1-ie^x}=\arctan e^x.
$$
The second derivative of the functional is the quadratic form
\begin{equation*}
  D^2 S =
\sum_{(j,k)}\frac{1}{\cosh(\rho_k-\rho_j)}\,(d\rho_k-d\rho_j)^2 ,
\end{equation*}
where the sum is taken over pairs of neighboring circles, which implies the
convexity.
\end{proof}

Now the idea is to minimize $S(\rho)$ restricted to $\sum_j\rho_j=0$. The convexity
of the functional implies the existence and uniqueness of solutions of the
classical boundary valued problems: Dirichlet (with prescribed $\rho_j$ on the
boundary) and Neumann (with prescribed $\Phi_j$ on the boundary).

The Schwarz-P surface is a triply periodic minimal surface. It is the symmetric
case in a $2$-parameter family of minimal surfaces with $3$ different hole sizes
(only the ratios of the hole sizes matter), see~\cite{DHKW92}. The Gauss map is a
double cover of the sphere with $8$ branch points. The image of the curvature line
pattern under the Gauss map is shown schematically in Fig.~\ref{f.recipe} {\em(top
left)}, thin lines. It is a refined cube. More generally, one may consider three
different numbers $m$, $n$, and $k$ of slices in the three directions. The $8$
corners of the cube correspond to the branch points of the Gauss map. Hence, not
$3$ but $6$ edges are incident with each corner vertex.  The corner vertices are
assigned the label $\cir$. We assume that the numbers $m$, $n$, and $k$ are even,
so that the vertices of the quad graph may be labelled `$\cir$', `$\sph$', and
`$\bullet$' consistently (see Section~\ref{s.isothermic}).

\begin{figure}[htb]%
\centering%
\includegraphics[width=0.8\textwidth]{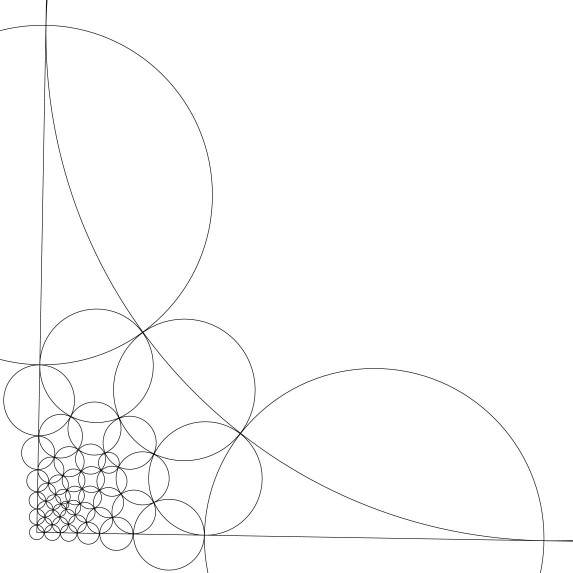}%
\caption{A piece of the circle pattern for a Schwarz-P surface after
  stereographic projection to the plane.}%
\label{f.euclid_schwarz_pattern}%
\end{figure}

We will take advantage of the symmetry of the surface and
construct only a piece of the corresponding circle pattern. Indeed
the combinatorics of the quad-graph has the symmetry of a
rectangular papallelepipedon. We construct an orthogonal circle
pattern with the symmetry of the rectangular papallelepipedon
eliminating the M\"obius ambiguity of Theorem~\ref{t.orthoKoebe}.
Consider one fourth of the sphere bounded by two orthogonal great
circles connecting the north and the south poles of the sphere.
There are two distinguished vertices (corners of the cube) in this
piece. Mapping the north pole to the origin and the south pole to
infinity by stereographic projection we obtain a Neumann boundary
value problem for orthogonal circle patterns in the plane. The
symmetry great circles become two orthogonal straight lines. The
solution of this problem is shown in
Fig.~\ref{f.euclid_schwarz_pattern}. The Neumann boundary data are
$\Phi=\pi/2$ for the lower left and upper right boundary circles
and $\Phi=\pi$ for all other boundary circles (along the symmetry
lines).

Now map this circle pattern to the sphere by the same
stereographic projection. One half of the spherical circle pattern
obtained (above the equator) is shown in Fig.~\ref{f.recipe}
{\em(top right)}. This is one eighth of the complete spherical
pattern. Now lift the circle pattern to the branched cover,
construct the Koebe polyhedron and dualize it to obtain the
Schwarz-P surface; see Fig.~\ref{f.recipe} {\em(bottom row)}. A
translational fundamental piece of the surface is shown in
Fig.~\ref{f.schwarz}.

We summarize these results in a theorem.

\begin{theorem}
  Given three even positive integers $m$, $n$, $k$, there exists a
  corresponding unique (unsymmetric) S-isothermic Schwarz-P surface.
\end{theorem}

Surfaces with the same ratios $m:n:k$ are different discretizations of the same
continuous Schwarz-P surface.  The cases with $m=n=k$ correspond to the symmetric
Schwarz-P surface.

Further examples of discrete minimal surfaces can be found in
\cite{BobenkoHoffmannSpringborn} and \cite{Buecking}.

Orthogonal circle patterns on the sphere are treated as the Gauss
map of the discrete minimal surface and are central in this
theory. Although circle patterns on the plane and on the sphere
differ just by the stereographic projection some geometric
properties of the Gauss map can get lost when represented in the
plane. Moreover to produce branched circle patterns in the sphere
it is important to be able to work with circle patterns directly
on the sphere. A variational method which works directly on the
sphere was suggested in
\cite{Springborn,BobenkoHoffmannSpringborn}. This variational
principle for spherical circle patterns is analogous to the
variational principles for Euclidean and hyperbolic patterns
presented in \cite{BobenkoSpringborn-Variational}. Unlike the
Euclidean and hyperbolic cases the spherical functional is not
convex, which makes difficult to use it in the theory. However the
spherical functional has proved to be amazingly powerful for
practical computation of spherical circle patterns (see
\cite{BobenkoHoffmannSpringborn} for detail). In particular, it
can be used to produce branched circle patterns in the sphere.

Numerous examples of discrete minimal surfaces are constructed with the help of the
spherical functional in the contribution of B\"ucking \cite{Buecking} to this
volume.

\section{Discrete Conformal Surfaces and Circle Patterns}
\label{s.conformal}

Conformal immersions
\begin{eqnarray*}
f:{\mathbb R}^2\supset D&\to& {\mathbb R}^3\\
 (x,y)  &\mapsto& f(x,y)
\end{eqnarray*}
are characterized by the properties
$$
\label{eq:isothermic} \|f_x\|=\|f_y\|, \quad f_x\bot f_y.
$$
Any surface can be conformally parametrized, and a conformal
parametrization induces a complex structure on $D$ in which
$z=x+iy$ is a local complex coordinate. The development of a
theory of discrete conformal meshes and discrete Riemann surfaces
is one of the popular topics in discrete differential geometry.
Due to their almost square quadrilateral faces and general
applicability, discrete conformal parametrizations are important
in computer graphics, in particular for texture mapping. Recent
progress in this area could be a topic of another paper; it lies
beyond the scope of this survey. In this section we mention
shortly only the methods based on circle patterns.

Conformal mappings can be characterized as mapping infinitesimal
circles to infinitesimal circles. The idea to replace
infinitesimal circles with finite circles is quite natural and
leads to circle packings and circle patterns (see also
Section~\ref{s.isothermic}). The corresponding theory for
conformal mappings to the plane is well developed (see the recent
book of Stephenson \cite{Stephenson}). The discrete conformal
mappings are constructed using a version of Koebe's theorem or a
variational principle which imply the corresponding existence and
uniqueness statements as well as a numerical method for
construction (see Section~\ref{s.minimal}). It is proven that a
conformal mapping can be approximated by a sequence of
increasingly fine, regular circle packings or patterns
\cite{Schramm-SG,HeSchramm}.

Several attempts have been made to generalize this theory for discrete conformal
parametrizations of surfaces.

The simplest natural idea is to ignore the geometry and take only
the combinatorial data of the simplicial surface
\cite{Stephenson}. Due to Koebe's theorem there exists essentially
unique circle packing representing this combinatorics. One treats
it as a discrete conformal mapping of the surface. This method has
been successfully applied by Hurdal et al. \cite{HurdalBSSRSR} for
flat mapping of the human cerebellum. However a serious
disadvantage is that the results depend only on the combinatorics
and not on the geometry of the original mesh.

An extension of Stephenson's circle packing scheme which takes the geometry into
account is due to Bowers and Hurdal \cite{BowersHurdal}. They treat circle patterns
with non-intersecting circles corresponding to vertices of the original mesh. The
geometric data in this case are the so called {\em inversive distances} of pairs of
neighboring circles, which can be treated as imaginary intersection angles of
circles. The idea is to get a discrete conformal mapping as a circle pattern in the
plane with the inversive distances coinciding with the inversive distances of small
spheres in space. The latter are centered at the vertices of the original meshes.
The disadvantage of this method is that there are almost no theoretical results
regarding the existence and uniqueness of inversive distance patterns.

Kharevich, Springborn and Schr\"oder suggested another way
\cite{KharevichSS} to handle the geometric data. They consider the
circumcircles of the faces of a simplicial surface and take their
intersection angles $\theta$. The circumcircles are taken
intrinsically, i.e. they are round circles with respect to the
surface metric. The latter is a flat metric with conical
singularities at vertices of the mesh. The intersection angles
$\theta$ are the geometric data to be respected, i.e. ideally for
a discrete conformal mapping one wishes a circle pattern in the
plane with the same intersection angles $\theta$. An advantage of
this method is that similarly to the circle packing method it is
based on a solid theoretical background --- the variational
principle for patterns of intersecting circles
\cite{BobenkoSpringborn-Variational}. However to get a circle
pattern flat one has to change the intersection angles
$\theta\mapsto\tilde\theta$ and it seems that there is no
geometric way to do this (in \cite{KharevichSS} the angles
$\tilde\theta$ of the circle pattern in the plane are defined as
minimizing the sum of differences squared $(\theta -
\tilde\theta)^2 $). A solution to this problem could possibly be
achieved by a method based on Delaunay triangulations of circle
patterns with disjoint circles. The corresponding variational
principle has been found recently in
\cite{Springborn-DelaunayCircle}.

It seems that discrete conformal surface parametrizations are at
the beginning of a promising development. Although now only some
basic ideas about discrete conformal surface parametrizations have
been clarified and no approximation results are known, there is
good chance for a fundamental theory with practical applications
in this field.


\bibliographystyle{amsplain}
\bibliography{BobenkoSC}

\end{document}